\documentclass[10pt]{article}
  \usepackage{epsfig,amsmath,amssymb,amsthm,cancel,bm,floatflt}
\usepackage{graphicx}
\usepackage{setspace} 
\usepackage{array}
\usepackage{pstricks,pstricks-add,pst-math,pst-xkey}

\newtheorem{theorem}{Theorem}

\newtheorem{corollary}[theorem]{Corollary}

\newtheorem{property}[theorem]{Property}

\newcommand{\real}{I\!\!R}

\newcommand{\btwobytwo}{\left [ \begin{array}{rr}}
\newcommand{\etwobytwo}{ \end{array} \right ]}
\newcommand{\bthrbythr}{\left [ \begin{array}{rrr}}
\newcommand{\ethrbythr}{ \end{array} \right ]}

\title{Intersections of hyperplanes\\ with conic sections in $\real^n$}

\author{P. M. Dearing \\ 
Department of Mathematical Sciences, Clemson University, \\
Clemson, SC 29632 \\
E-mail: pmdrn@clemson.edu}

\begin{document}
\maketitle

\begin{abstract}
Closed form expressions are given for computing the parameters and vectors that identify and define the $n-1$ dimensional conic section
that results from the intersection of a hyperplane with each of the  $n$-dimensional conic section: cone, hyperboloid of two sheets, ellipsoid or paraboloid.  The conic sections are assumed to be symmetric about their major axis, but may have any orientation and center.  A class of hyperboloids are identified with the property that the parameters and vectors of the intersection of all hyperboloids in a subset of the class can be computed efficiently.
\end{abstract}

\section{Introduction}

This paper considers the intersection of a hyperplane in $\real^n$ with each of the $n$-dimensional conic sections: cone, hyperboloid of two sheets, ellipsoid, or paraboloid. Each conic section is assumed to be symmetric about its major axis, but may have any orientation and center.  Closed form expressions are given for identifying the resulting intersection as an $n-1$ dimensional hyperboloid, ellipsoid or paraboloid, and for computing the parameters and vectors that define the resulting conic section. 

Also considered is a class of $n$-dimensional hyperboloids defined by a finite set $\mathcal{P}$ of points in $\real^n$, so that for each point $\mathbf{p}_i \in \mathcal{P}$, there is a corresponding non-negative number $r_i$, and the Euclidean ball $[\mathbf{p}_i, r_i]$, with center $\mathbf{p}_i$ and radius $r_i$.  Conditions are given so that each pair of balls $[\mathbf{p}_j, r_j]$, $[\mathbf{p}_k, r_k]$, with $r_j > r_k$, determines a hyperboloid with focal points $\mathbf{p}_j$ and $\mathbf{p}_k$, and constant $r_j - r_k$.  

This class of hyperboloids has the property that for a subset $\mathcal{S} \subseteq \mathcal{P}$ of size $k$, the intersection of hyperboloids corresponding to all pairs of points is $\mathcal{S}$ can be computed as the intersection of a sequence of $k-1$ hyperboloids with a common  focal point.

For each subset $\mathcal{T} \subseteq \mathcal{S}$ of three points with distinct radii, and the three hyperboloids corresponding to each pair of points in $\mathcal{T}$, a hyperplane is constructed with the property that the intersection of any two of the hyperboloids equals the intersection of the hyperplane with either of the two hyperboloids.
This property is then used to show that the intersection of a sequence of $k-1$ hyperboloids is equal to the intersection of one of the $k-1$ hyperboloids with $k-2$ hyperplanes.  The resulting intersection yields a conic section of dimension $n-k+2$,  along with the vectors and parameters that define the conic section.

The results presented here were motivated by two applications.  One is the optimization problem of finding the $n$-dimensional ball of minimum radius  that contains a finite set of $n$-dimensional balls.  In  \cite{Dearingmincov}, primal and dual algorithms are constructed for this problem based on finding the intersection of a sequence of hyperboloids.

In \cite{Dearingtdoa}, an alternative solution approach is presented to the problem of  locating the source (e.g. cell phone) of an electronic signal  that is transmitted to a set of receivers.  The problem data consists of the receiver locations and differences in the time of arrival of the signal at pairs of receivers.   
In \cite{leva}, Leva presents a solution approach to this problem based on the intersection of hyperboloids in $\real^3$.  Reference \cite{Dearingtdoa} expands the approach in \cite{leva} to the intersections of all conic sections in $\real^n$, and to a determination of 
when there is a unique or alternate solution. 

The following sections present known definitions and characterizations for each of the $n$-dimensional conic sections \cite{Brannan}, \cite{Glaeser}. Then a characterization of conic sections is presented in terms of a quadratic form that is  used to analyze the intersections of a hyperplane with each conic section.  A procedure to determine the intersection of a sequence of hyperboloids is presented, and properties of a class of hyperboloids is presented.
 
\section{Hyperboloids}
\pagenumbering{arabic} \label{sec:intro}

An $n$-dimensional \textbf{hyperboloid} of two sheets, symmetric about its major axis, is the set $H$ of all points $\mathbf{x} \in \real^n$ such that the absolute difference in the distance from $\mathbf{x}$ to two given  points $\mathbf{p}_1, \mathbf{p}_2$, equals a positive constant $2a$. That is,
\begin{equation}
H = \{ \mathbf{x}: | \parallel \mathbf{p}_{2} - \mathbf{x} \parallel - \parallel \mathbf{p}_{1} - \mathbf{x} \parallel | =  2a  \}. \label{H}
\end{equation}
The \textbf{sheet} of $H$ closest to $\mathbf{p}_1$ is the set
\begin{equation}
H_{1} = \{ \mathbf{x}: \parallel \mathbf{p}_{2} - \mathbf{x} \parallel - \parallel \mathbf{p}_{1} - \mathbf{x} \parallel  =  2a  \}, \label{H1}
\end{equation}
and the sheet closest to $\mathbf{p}_2$ is the set
\begin{equation}
H_{2} = \{ \mathbf{x}: \parallel \mathbf{p}_{2} - \mathbf{x} \parallel - \parallel \mathbf{p}_{1} - \mathbf{x} \parallel  =  -2a  \}. \label{H2}
\end{equation}
Observe that  $H = H_1 \cup H_2$.

A hyperboloid $H$ is specified by the following vectors and parameters, all of which are determined by the points $\mathbf{p}_1$ and $\mathbf{p}_2$,
and the positive constant $a$. 
The \textbf{focal points} of  $H$ are the points $\mathbf{p}_1$ and $\mathbf{p}_2$. 
The \textbf{center} of $H$ is the mid-point of the line segment  between the focal points, $\mathbf{c} =\frac{1}{2}(\mathbf{p}_{1} + \mathbf{p}_{2})$,   and  the parameter $c =\frac{1}{2}\!\!\!\parallel\!\mathbf{p}_{1} - \mathbf{p}_{2}\!\!\parallel$ is the distance from the center to either focal point.  
The unit vector $\mathbf{v} = \frac{\mathbf{p}_{1} - \mathbf{p}_{2}}{\parallel \mathbf{p}_{1} - \mathbf{p}_{2} \parallel }$ 
 from $\mathbf{p}_{2}$ to  $ \mathbf{p}_{1}$ is called the \textbf{axis vector}, and is parallel to the \textbf{major axis},
 which is the line through  $ \mathbf{p}_{1}$ and $\mathbf{p}_{2}$. 
 The \textbf{vertex} of the sheet $H_{1}$  is the point $\mathbf{a}_{1} = \mathbf{c} + a\mathbf{v}$, and is the point of intersection between  the major axis and $H_{1}$.
 The vertex of $H_2$ is the point $\mathbf{a}_{2} = \mathbf{c} - a\mathbf{v}$, and is point of intersection between the major axis and $H_2$.  The \textbf{eccentricity} specifies the shape of $H$ and is given by  $\epsilon = \frac{c}{a}$.

The triangle inequality implies 
$ 2a = |\parallel\mathbf{p}_{2} - \mathbf{x}\parallel - \parallel\mathbf{p}_{1} - \mathbf{x}\parallel| \leq  \parallel\mathbf{p}_{2} -   \mathbf{p}_{1}\parallel = 2c $,
so that $a \leq c$.
If $a = c$,
then each sheet of the hyperboloid $H$ consists of a ray along the major axis from the  points $\mathbf{p}_1$ and $\mathbf{p}_2$ respectively, and is a degenerate hyperboloid.  That is, $H_1 = \{ \mathbf{x} = \mathbf{p}_1 + \alpha\mathbf{v}, \alpha \geq 0 \}$, and $H_2 = \{ \mathbf{x} = \mathbf{p}_2 - \alpha\mathbf{v}, \alpha \geq 0 \}$. For a hyperboloid $H$ it is assumed that $a < c$, or equivalently, that $\epsilon > 1$. Since $c > 0$ and fixed, and $0 < a < c$, then $1 < \epsilon < \infty$,

The \textbf{directrix} of the sheet $H_{1}$ is the hyperplane with normal vector $\mathbf{v}$  
containing  the point $\mathbf{d}_{1}$,
where $\mathbf{d}_{1} = \mathbf{c} + d \mathbf{v}$,  and 
$d = \frac{a^2}{c} = \frac{2a^2}{ \parallel  \mathbf{p}_{1}-\mathbf{p}_{2} \parallel}$.
The directrix of $H_1$ is defined by  $\mathbf{v} \mathbf{x}  = \mathbf{v}\mathbf{d}_{1} 
= \mathbf{v}\mathbf{c} + d
= \frac{ \parallel \mathbf{p}_{1} \parallel^2 -  \parallel \mathbf{p}_{2} \parallel^2 +  4a^2} {2  \parallel  \mathbf{p}_{1}-\mathbf{p}_{2} \parallel}$. 
 For the sheet $H_2$,  $\mathbf{d}_{2} = \mathbf{c} - d \mathbf{v}$, and the directrix of $H_2$  is defined by  $\mathbf{v} \mathbf{x}  = \mathbf{v}\mathbf{d}_{2} 
= \mathbf{v}\mathbf{c} - d = \frac{ \parallel \mathbf{p}_{1} \parallel^2 -  \parallel \mathbf{p}_{2} \parallel^2 -  4a^2} {2  \parallel  \mathbf{p}_{1}-\mathbf{p}_{2} \parallel}$. 
The first property states a well known equivalent expression for each sheet $H_1$ and $H_2$ in terms of its directrix. 
\begin{property}  The sheet $H_1 = H^*_1$ and the sheet $H_2 = H^*_2$ where 
\begin{alignat}{3}
&H^*_1 = &\{ \mathbf{x}: \parallel \mathbf{p}_{1} - \mathbf{x} \parallel = \epsilon( \mathbf{v} \mathbf{x} - \mathbf{v} \mathbf{d}_1 ) \},  \label{HH1} \\
&H^*_2 = &\{ \mathbf{x}: \parallel \mathbf{p}_{1} - \mathbf{x} \parallel = \epsilon( \mathbf{v} \mathbf{d}_2 - \mathbf{v} \mathbf{x} ) \}.  \label{HH2} 
\end{alignat}
\end{property}
\noindent \textbf{Proof:} The proof expands \eqref{H1} and \eqref{H2} and substitutes the definitions of parameters and vectors to obtain the results. \hfill $\Box$

For any unit vector $\mathbf{u}$ orthogonal to the axis vector $\mathbf{v}$ of a hyperboloid $H$ with center $\mathbf{c}$, the space curve
$\{ \mathbf{x}_1(\alpha) = \mathbf{c} + a \sec(\alpha) \mathbf{v} + b \tan (\alpha) \mathbf{u}, - \pi/2 < \alpha < \pi/2 \}$, where $b = \sqrt{c^2 - a^2}$ 
gives a parametric representation of one sheet of a hyperbola in the 
two dimensional affine space aff$(\mathbf{v},\mathbf{u},\mathbf{c})$.  The next property shows that $\mathbf{x}_1(\alpha)$ is a subset of the sheet $H_1$.  An analogous result holds for the sheet $H_2$.

\begin{property}
Given a hyperboloid $H$ in $\real^n$ with focal points $\mathbf{p}_1$ and $\mathbf{p}_2$, center $\mathbf{c}$, axis vector $\mathbf{v}$, eccentricity $\epsilon$, and sheets $H_1$ and $H_2$,  
if $\mathbf{u}$ is a unit vector orthogonal to the axis vector $\mathbf{v}$,
then  $H \cap$aff$(\mathbf{v},\mathbf{u},\mathbf{c})$ is a two dimensional hyperbola with the same focal points, center, axis, and eccentricity as $H$.
Furthermore, the two-dimensional sheets of $H \cap$aff$(\mathbf{v},\mathbf{u},\mathbf{c})$ are given by
 $H_1 \cap $aff$(\mathbf{v},\mathbf{u},\mathbf{c}) = \{ \mathbf{x}_1(\alpha) = \mathbf{c} + a \sec(\alpha) \mathbf{v} + b \tan (\alpha) \mathbf{u}, - \pi/2 < \alpha < \pi/2 \}$, where $b = \sqrt{c^2 - a^2}$, and 
$H_2\cap$aff$(\mathbf{v},\mathbf{u},\mathbf{c}) = \{ \mathbf{x}_2(\alpha) = \mathbf{c} - a \sec(\alpha) \mathbf{v} + b \tan (\alpha) \mathbf{u}, - \pi/2 < \alpha < \pi/2 \} = \{ \mathbf{x}_2(\alpha) = \mathbf{c} + a \sec(\alpha) \mathbf{v} + b \tan (\alpha) \mathbf{u},  \pi/2 < \alpha < 3\pi/2 \}$.
 \end{property}
 \noindent \textbf{Proof:} The proof  shows that $\mathbf{x}_1(\alpha)$ satisfies expression \eqref{HH1} (and $\mathbf{x}_2(\alpha)$ satisfies expression \eqref{HH2}) by expanding \eqref{HH1} and \eqref{HH2} and substituting definitions of parameters and vectors to obtain the results. \hfill $\Box$


\section{Ellipsoids}

An $n$-dimensional \textbf{ellipsoid},  symmetric about its major axis, is the set $E$ of all points $\mathbf{x} \in \real^n$ such that the sum of the distances from $\mathbf{x}$ to two given  points $\mathbf{p}_1$ and $\mathbf{p}_2$, equals a positive constant $2a$. That is,
\begin{equation}
E = \{ \mathbf{x}: \parallel \mathbf{p}_{2} - \mathbf{x} \parallel + \parallel \mathbf{p}_{1} - \mathbf{x} \parallel =  2a  \}. \label{E}
\end{equation}

An ellipsoid $E$ is specified by the same vectors and parameters that specify a hyperboloid, all of which are  determined by the focal points $\mathbf{p}_1$ and $\mathbf{p}_2$ 
and the positive constant $a$.
That is, the  axis vector $\mathbf{v}$, the center point $\mathbf{c}$, the parameter $c$, the vertices $\mathbf{a}_1$ and 
$\mathbf{a}_2$, the eccentricity $\epsilon =  \frac{c}{a}$, and the directrix $\mathbf{d}_1$ and $\mathbf{d}_2$, each have the same definition for an ellipsoid as for a hyperboloid, except that $c \leq a$, as shown next.

Using the triangle inequality, 
$ 2c = \parallel\mathbf{p}_{1} -   \mathbf{p}_{2}\parallel \leq \parallel\mathbf{p}_{1} - \mathbf{x}\parallel + \parallel\mathbf{p}_{2} - \mathbf{x}\parallel = 2a $,
so that $c \leq a$.
If $a = c$, $E$  is the line segment between $\mathbf{p}_1$ and $\mathbf{p}_2$, and is a degenerate ellipsoid.   
 Thus  for an ellipsoid  $E$ it is assumed that $c < a$.  Since  $0 \leq c < a$, then $0 \leq \epsilon < 1$.
 
 The first property states two well known equivalent expressions for an ellipsoid  $E$ in terms of the  directrix. 
\begin{property}   $E = E_1$ and  $E = E_2$ where.\begin{alignat}{3}
&E_1 = &\{ \mathbf{x}: \parallel \mathbf{p}_{1} - \mathbf{x} \parallel = \epsilon( \mathbf{v} \mathbf{d}_1 - \mathbf{v} \mathbf{x} ) \}, \label{E1} \\
&E_2 = &\{ \mathbf{x}: \parallel \mathbf{p}_{2} - \mathbf{x} \parallel = \epsilon( \mathbf{v} \mathbf{x} - \mathbf{v} \mathbf{d}_2 ) \}.  \label{E2} 
\end{alignat}
\end{property}
\noindent \textbf{Proof:} To show $E \subseteq E_1$, write  \eqref{E} as 
$\parallel\mathbf{p}_1 - \mathbf{x} \parallel - 2a = -\parallel \mathbf{p}_2 - \mathbf{x} \parallel$, square both sides and substitue parameters and vectors of $H$.
A reverse argument shows that $E_1 \subseteq E$.
 An analogous argument shows $E = E_2$. \hfill $\Box$

\begin{property}
Given an ellipsoid $E$ with center $\mathbf{c}$, axis vector $\mathbf{v}$ and eccentricity $\epsilon$, if $\mathbf{u}$ is a unit vector orthogonal to the axis vector $\mathbf{v}$, then $E\cap$aff$(\mathbf{u}, \mathbf{v}, \mathbf{c}) = \{ \mathbf{x}(\alpha) = \mathbf{c} + a \cos(\alpha) \mathbf{v} + b \sin (\alpha) \mathbf{u}, 0 < \alpha < 2\pi \}$, where $b = \sqrt{a^2 - c^2}$,  is a two-dimensional ellipse with the same center, axis vector and eccentricity as $E$.  \end{property}
 \noindent \textbf{Proof:} The proof shows that $\mathbf{x}(\alpha)$ satisfies expressions \eqref{E1} (and \eqref{E2}) by expanding 
 \eqref{E1} and \eqref{E2} and substituting definitions of parameters and vectors to obtain the results. \hfill $\Box$


\section{Quadratic form representation of hyperboloids,\\ ellipsoids and cones}

The next property gives an equivalent representation for a hyperboloid or an ellipsoid in terms of a quadratic form. 
An equivalent form of this representation for hyperboloids in $\real^3$ is reported in \cite{leva}
\begin{property} Given the focal points $\mathbf{p}_1$ and $\mathbf{p}_2$ and a positive constant $a$, with corresponding axis vector $\mathbf{v}$, center point $\mathbf{c}$, and eccentricity $\epsilon = c/a$, let $H$ be the hyperboloid determined by these vectors and parameters if $c > a$, and let $E$ be the ellipsoid determined by these vectors and parameters if $c< a$. Let $Q$ be the set defined by 
\begin{equation}
Q = \{  \mathbf{x}:  ( \mathbf{x} - \mathbf{c})^T [ I - \epsilon^2 \mathbf{v} \mathbf{v}^T] ( \mathbf{x} - \mathbf{c})  = a^2 - c^2 \}.
\label{HQ1} \end{equation}

Then $Q = H$ if $c > a$, and $Q = E$ if $c < a$. 
\end{property}
\noindent \textbf{Proof:} Assume that $c > a$. To prove that $H_1 \subseteq Q$, square both sides of expression \eqref{HH1} and substitute the parameters and vectors of the  hyperboloid $H$ to obtain expression \eqref{HQ1}.  Observe that $( \mathbf{x} - \mathbf{c})^T [I -  \epsilon^2 \mathbf{v} \mathbf{v}^T] ( \mathbf{x} - \mathbf{c}) =  [ \mathbf{x} - \mathbf{c}]^2 -  [ \epsilon \mathbf{v}( \mathbf{x} - \mathbf{c})]^2$.
The proof that $H_2 \subseteq Q$ is analogous, starting with expression \eqref{HH2} and squaring both sides.  Thus  $H_1 \cup H_2 = H \subseteq Q$.

Assume that $c < a$. To prove that $E_1 \subseteq Q$, expand expression \eqref{E1}, square both sides, and substitute the parameters and vectors of the  ellipsoid $E$ to obtain expression \eqref{HQ1}.  A similar argument starting with \eqref{E2} shows that $E_2 \subseteq Q$.

To show reverse set inclusion, assume that $\mathbf{x} \in Q$ and apply the argument above in reverse, introducing either $\mathbf{p}_1$ or $\mathbf{p}_2$. 
One case leads to $\mathbf{x} \in E_1 = E$ and $c < a$.  The other case leads to  $\mathbf{x} \in H$ and $c > a$.  
Thus $Q = E$ with $c < a$, or $Q = H$ with $c > a$.\hfill $\Box$

Observe that  $[I - \epsilon^2 \mathbf{v} \mathbf{v}^T]$ is similar to a Householder matrix. Direct computation shows that 
the matrix  in \eqref{HQ1} has eigenvalue $1-\epsilon^2$ of multiplicity one with  $\mathbf{v}$
as the corresponding eigenvector,   
and that $1$ is an eigenvalue of multiplicity $n-1$, with corresponding eigenvectors orthogonal to  $\mathbf{v}$ and mutually orthogonal.   
Thus  $Q$ is symmetric about its major axis. Observe that if $Q$ is an ellipsoid, so that $\epsilon < 1$, then the matrix $ [ I - \epsilon^2 \mathbf{v} \mathbf{v}^T ]$ is positive definite, and the right hand side of \eqref{HQ1} is positive, which is consistent with the the representation of an ellipsoid in terms of a quadratic form with positive definite matrix.

An $n$-dimensional right circular cone, denoted by $C$, with center $\mathbf{c}$, axis vector $\mathbf{v}$, and eccentricity $\epsilon$, may also be expressed in terms of the  quadratic form $Q$ with $\epsilon > 1$, but with a right hand side value of zero.  That is,
\begin{equation}
C = \{ \mathbf{x} : (\mathbf{x} - \mathbf{c})^T[I - \epsilon^2 \mathbf{v} \mathbf{v}^T ] (\mathbf{x} - \mathbf{c}) = 0 \}.
\label{C} \end{equation}
Note that $C$ is a cone since if $\mathbf{x} - \mathbf{c} \in C$, then $\lambda(\mathbf{x} - \mathbf{c}) \in C$ for $\lambda \geq 0$.

The cone $C$ has two "sheets" denoted by $C_1$ and $C_2$, where $C_1 = \{ \mathbf{x} : \parallel \mathbf{x} - \mathbf{c} \parallel = \epsilon \mathbf{v}(\mathbf{x} - \mathbf{c}) \}$,  and is the subset of $C$ that is closest to the focal point $\mathbf{p}_1$.  The sheet $C_2 = \{ \mathbf{x} : \parallel \mathbf{x} - \mathbf{c} \parallel = \epsilon \mathbf{v}(\mathbf{c} - \mathbf{x}) \}$ and is the subset of $C$ closest to the focal point $\mathbf{p}_2$.  Observe that $\parallel \mathbf{x} - \mathbf{c} \parallel/(\mathbf{x} - \mathbf{c})\mathbf{v} = c/a = \epsilon$.   

\begin{property}  $C = C_1 \cup C_2$.
\end{property}
 \noindent \textbf{Proof:} If $\mathbf{x} \in C_1$, then 
  $\parallel \mathbf{x} - \mathbf{c} \parallel^2 = \epsilon^2 [\mathbf{v}( \mathbf{x} - \mathbf{c})]^2$, and $(\mathbf{x} - \mathbf{c})^T[I -  \epsilon^2 \mathbf{v} \mathbf{v}^T] (\mathbf{x} - \mathbf{c}) = 0$, so that $C_1 \subseteq C$.  A similar argument shows that $C_2 \subseteq C$. Applying the argument in reverse shows that if $\mathbf{x} \in C$, then 
 $\parallel \mathbf{x} - \mathbf{c} \parallel^2 = \epsilon^2 [\mathbf{v}( \mathbf{x} - \mathbf{c})]^2$ so that $\parallel \mathbf{x} - \mathbf{c} \parallel = \pm \epsilon [\mathbf{v}( \mathbf{x} - \mathbf{c})]$, and either $\mathbf{x} \in C_1$, or $\mathbf{x} \in C_2$.  \hfill  $\Box$
 
 For any point $\mathbf{x} \in C$, let $\gamma$ be the angle between the vector $\mathbf{x} - \mathbf{c}$ and the axis vector.  Then  $\sec(\gamma) = \frac{c}{a} = \epsilon$.
The next Property and its proof are analogous to Property 2.

\begin{property}
Given a cone $C$ in $\real^n$ with center $\mathbf{c}$, axis vector $\mathbf{v}$, and sheets $C_1$ and $C_2$,  
if $\mathbf{u}$ is a unit vector orthogonal to  $\mathbf{v}$,
then $C\cap$aff$(\mathbf{u}, \mathbf{v}, \mathbf{c})$ is a two-dimensional cone with the same center and axis vector as $C$.  Furthermore, 
the two dimensional sheets of $C\cap$aff$(\mathbf{u}, \mathbf{v}, \mathbf{c})$ are given by
 $C_1\cap $aff$(\mathbf{u}, \mathbf{v}, \mathbf{c}) =   \{ \mathbf{x}'_1(\beta) = \mathbf{c} + a|\beta| \mathbf{v} + b \beta \mathbf{u}, - \infty < \beta < \infty \}$, where $b = \sqrt{c^2 - a^2}$, and 
  $C_2\cap$aff$(\mathbf{u}, \mathbf{v}, \mathbf{c}) =   \{ \mathbf{x}'_2(\beta) = \mathbf{c} - a|\beta| \mathbf{v} + b \beta \mathbf{u}, - \infty < \beta < \infty \}$.
 \end{property}
  \noindent \textbf{Proof:}  Substituting $\mathbf{x}'_1(\beta)$ into the left hand side of the expression for $C_1$ and squaring gives the right hand side of the expression for $C_1$.
   The same approach shows $\mathbf{x}'_2(\beta)$ satisfies the expression for $C_2$.  \hfill $\Box$
 
\begin{property} If a cone $C$ and a hyperboloid $H$ have the same axis vector $\mathbf{v}$, center $\mathbf{c}$, and eccentricity $\epsilon$, then $C$ is the asymptotic approximation of $H$.
\end{property}
\noindent \textbf{Proof}: Choose any  unit vector $\mathbf{u}$ orthogonal to $\mathbf{v}$, and consider the sheet of the two-dimensional hyperbola $\mathbf{x}_1(\alpha)$ and the sheet of the two-dimensional cone $\mathbf{x}_1'(\beta)$.
In $\mathbf{x}_1'(\beta)$, 
substitute the parameter $\tan(\alpha)$, for $- \pi/2 < \alpha < \pi/2$, in place of $\beta$, for $- \infty < \beta < \infty $.
Then $\lim_{\alpha \rightarrow  |\pi/2|} \parallel \mathbf{x}'_1(\alpha) - \mathbf{x}_1(\alpha)\parallel 
= a \parallel\mathbf{v} \parallel \lim_{\alpha \rightarrow |\pi/2|}  (|\tan(\alpha)| - \sec(\alpha)) = 0$,
which shows that the cone $\mathbf{x}_1'(\alpha)$ is the asymptotic approximation to the hyperbola $\mathbf{x}_1(\alpha)$.
A similar analysis shows that $\mathbf{x}_2'(\alpha)$ is the asymptotic approximation to the hyperbola $\mathbf{x}_2(\alpha)$.
Since these results hold for any  unit vector $\mathbf{u}$ orthogonal to $\mathbf{v}$, the cone 
$C$ is the asymptotic approximation to the hyperboloid $H$. \hfill $\Box$


\section{Parabaloids}

An $n$-dimensional \textbf{paraboloid}, symmetric about its major axis, is the set $P$ of all points $\mathbf{x} \in \real^n$ such that the distance from $\mathbf{x}$ to a given point $\mathbf{p}_1$ on the major axis, equals  the distance from $\mathbf{x}$ to the hyperplane that is orthogonal to the major axis and contains a point $\mathbf{p}_2$ on the major axis.
A paraboloid is specified by the two points $\mathbf{p}_1$ and  $\mathbf{p}_2$ only.
The point $\mathbf{p}_1$ is the focal point of the paraboloid.
The  major axis is the line through the points $\mathbf{p}_1$ and  $\mathbf{p}_2$.   The axis vector $\mathbf{v} = \frac{(\mathbf{p}_1 - \mathbf{p}_2)}{\parallel \mathbf{p}_1 - \mathbf{p}_2 \parallel}$ is the unit vector parallel to the major axis.  
A paraboloid $P$, defined by the points $\mathbf{p}_1$ and $\mathbf{p}_2$ in $\real^n$,
is the set  
\begin{equation}
P = \{ \mathbf{x}: \parallel \mathbf{p}_{1} - \mathbf{x} \parallel =  \mathbf{v}(\mathbf{x} - \mathbf{p_2}) \}. \label{P}
\end{equation}

The vertex of the paraboloid  is the the center $\mathbf{c}$, where  $\mathbf{c} = \frac{\mathbf{p}_1 + \mathbf{p}_2}{2}$, and is the intersection of the paraboloid with the major axis. 
The  parameter $c = \frac{\parallel \mathbf{p}_1 - \mathbf{p}_2 \parallel}{2}$. 
The \textbf{directrix} of a paraboloid $P$ is the hyperplane with normal vector $\mathbf{v}$  containing  the point $\mathbf{p}_{2}$.  Observe that $\mathbf{c} - \mathbf{p}_2 = \frac{\mathbf{p}_1 + \mathbf{p}_2 - \mathbf{p}_2 - \mathbf{p}_2}{2} = \frac{\mathbf{p}_1 - \mathbf{p}_2}{2}\frac{\parallel \mathbf{p}_1 - \mathbf{p}_2 \parallel}{\parallel \mathbf{p}_1 - \mathbf{p}_2 \parallel} = c \mathbf{v}$, and $\mathbf{c} - \mathbf{p}_1 = -c\mathbf{v}$.  
There is no parameter $a$ for a paraboloid, and all paraboloids have the same shape, so there is no parameter like the eccentricity for hyperboloids and ellipsoids.

For any unit vector $\mathbf{u}$ orthogonal to the axis vector $\mathbf{v}$ of a parabaloid $P$ with center $\mathbf{c}$, the space curve
$\{ \mathbf{x}(\alpha) = \mathbf{c} + c \alpha^2 \mathbf{v} + 2c \alpha \mathbf{u}, - \infty < \alpha < \infty \}$ 
gives a parametric representation of a parabola in the 
two dimensional affine space aff$(\mathbf{v},\mathbf{u},\mathbf{c})$.  The next property shows that $\mathbf{x}(\alpha)$ is a subset of the paraboloid  $P$. 

\begin{property}
Given a paraboloid $P$ in $\real^n$, with focal points $\mathbf{p}_1$ and $\mathbf{p}_2$  
if $\mathbf{u}$ is a unit vector orthogonal to the axis vector $\mathbf{v}$,
then the space curve
\begin{equation}
\{ \mathbf{x}(\alpha) = \mathbf{c} + c \alpha^2 \mathbf{v} + 2c \alpha \mathbf{u}, - \infty < \alpha < \infty \} \label{SE1}
 \end{equation}
 is a subset of the paraboloid $P$ and is a two dimensional parabola in the affine space aff$(\mathbf{v}, \mathbf{c}, \mathbf{u})$ with focal points 
 $\mathbf{p}_1$ and $\mathbf{p}_2$.
 \end{property}
 \noindent \textbf{Proof:} The proof shows that $\mathbf{x}(\alpha)$ satisfies expression \eqref{P} by expansion and substitution of \eqref{P}. \hfill $\Box$

Paraboloids may also be expressed  in terms of a quadratic form similar to \eqref{HQ1}. 
However, for a paraboloid, there is no eccentricity, and the right hand side is a linear expression of $\mathbf{x}$.
\begin{property}
The paraboloid $P$ has the equivalent expression.
\begin{equation}
P = \{  \mathbf{x}:  ( \mathbf{x} - \mathbf{c})^T [I - \mathbf{v} \mathbf{v}^T] ( \mathbf{x} - \mathbf{c})  = 4c\mathbf{v}(\mathbf{x} - \mathbf{c}) \}.
\label{PQ} \end{equation}
\end{property}
\noindent \textbf{Proof:} Expanding and squaring the left hand side of \eqref{P} yields
$\parallel \mathbf{p}_1 -\mathbf{x} \parallel^2 =  \parallel \mathbf{x} - \mathbf{c} + \mathbf{c} -\mathbf{p}_1 \parallel^2 =
\parallel  \mathbf{x} - \mathbf{c} -c\mathbf{v} \parallel^2
= \parallel  \mathbf{x} - \mathbf{c} \parallel^2 - 2c\mathbf{v}( \mathbf{x} - \mathbf{c}) + c^2 $.
Expanding and squaring the right hand side of \eqref{P} yields
$[\mathbf{v}(\mathbf{x} - \mathbf{p}_2)]^2 = [\mathbf{v}(\mathbf{x} - \mathbf{c} + \mathbf{c} -\mathbf{p}_2)]^2 
 = [\mathbf{v}(\mathbf{x} - \mathbf{c}) + \mathbf{v}(c \mathbf{v})]^2 \notag \\
 = [\mathbf{v}(\mathbf{x} - \mathbf{c})]^2 +2c\mathbf{v}(\mathbf{x} - \mathbf{c}) +c^2
=(\mathbf{x} - \mathbf{c})^T[\mathbf{v}^T\mathbf{v}](\mathbf{x} - \mathbf{c}) +2c\mathbf{v}(\mathbf{x} - \mathbf{c}) +c^2$.
Equating the two sides and re-arranging yields 
$(\mathbf{x} - \mathbf{c})^T[I - \mathbf{v}^T\mathbf{v}](\mathbf{x} - \mathbf{c}) = 4c\mathbf{v}(\mathbf{x} - \mathbf{c})$. \hfill $\Box$


\section{Intersections of a hyperplane with a cone or hyperboloid}

   From the classical studies of conic sections in $\real^3$, it is well known that if a  plane  and a cone  
  intersect at an appropriate angle, measured between the axis vector of the cone and the normal vector of the plane, the intersection is either a two dimensional hyperbola, ellipse, or parabola.
  These results extend to the intersection of a hyperplane with each conic section in $\real^n$, and are reported below.  For the intersection of a hyperplane with a hyperboloid or a cone, conditions are given for the resulting intersection to be a hyperboloid, an ellipsoid or a paraboloid of dimension $n-1$. For the intersection of a hyperplane and an ellipsoid, the resulting intersection is always an ellipsoid of dimension $n-1$.
 For the intersection of a hyperplane and paraboloid, conditions are given for the resulting intersection to be a paraboloid or an ellipsoid of dimension $n-1$. 
 
  For the intersection of a hyperplane with each of the conic sections in $\real^n$, expressions are given for identifying the resulting  hyperboloid, ellipsoid or paraboloid of dimension $n-1$, and for computing the vectors and  parameters that characterize it.

\begin{property}  Suppose $Q = \{ \mathbf{x}: (\mathbf{x} - \mathbf{c})^T[ I - \epsilon^2 \mathbf{v} \mathbf{v}^T ] (\mathbf{x} - \mathbf{c}) = a^2 -c^2  \}$ 
 is a hyperboloid in $\real^n$,  centered at $\mathbf{c} = (c_1, \dots, c_n)^T$, 
with  axis  vector $\mathbf{v}$ of unit length,  eccentricity $\epsilon > 1$, and 
  parameters $a$ and $c$, and suppose 
  $H\!P = \{ \mathbf{x}: \mathbf{h}( \mathbf{x} -\mathbf{c})=\hat{h} \}$,
  is a hyperplane  with  $\| \mathbf{h} \| =1$.
  Let $\rho = \sqrt{1-(\mathbf{h} \mathbf{v})^2}$.  
 Then $Q \cap H\!P$ is a hyperboloid of dimension $n-1$ iff $\epsilon\rho > 1$, 
 or an ellipsoid of dimension $n-1$ iff $\epsilon\rho < 1$ and $\hat{h}^2 \geq a^2(1-\rho^2\epsilon^2)$, or a paraboloid of dimension $n-1$ iff  $\epsilon \rho = 1$.
 \end{property}
 \noindent \textbf{Proof:} 
 If $\mathbf{h}$ and $\mathbf{v}$ are linearly dependent, then $H\!P$ is orthogonal to $\mathbf{v}$ and $Q \cap H\!P$ is a ball of dimension $n-1$ with center $\mathbf{c} +\hat{h}\mathbf{h}$.
 Assume that  $\mathbf{h} $ and $\mathbf{v}$ are linearly independent, and let sub$(\mathbf{v},\mathbf{h})$ denote the  2-dimensional subspace generated by  
$\mathbf{h}$ and $\mathbf{v}$.
Let $\{ \mathbf{g}_1, \ldots, \mathbf{g}_{n-1} \}$ be an orthonormal basis of the  
 null space of $\mathbf{h}$ in $\real^n$, and suppose that  $\mathbf{g}_1$  is chosen  so that  $\mathbf{g}_1\mathbf{v} >0$,
 and so that $\mathbf{g}_1$ lies in  $\text{sub}(\mathbf{v}, \mathbf{h})$. 
 That is,  $\mathbf{g}_1 =( \mathbf{v}-\hat{\mathbf{g}})/ \| \mathbf{v}-\hat{\mathbf{g}} \| $,
   where $\hat{\mathbf{g}} = (\mathbf{v} \mathbf{h}) \mathbf{h}$ is the projection of $\mathbf{v}$ onto $\mathbf{h}$.
   
 Let $T$ be the $n \times n$ orthonormal matrix with rows  $ \mathbf{g}^T_1,\ldots, \mathbf{g}^T_{n-1}, \mathbf{h}^T $.  
  Then  $T(\mathbf{h}) = \boldsymbol{\varepsilon}_n$, and 
  $T(\mathbf{v}) = \mathbf{v}' = \rho\boldsymbol{\varepsilon}_1 + \sigma \boldsymbol{\varepsilon}_n$, 
   where $\boldsymbol{\varepsilon}_i$ is the $i^{th}$ unit vector,  $\rho  = \mathbf{g}_1\mathbf{v} > 0$, $\sigma=\mathbf{h} \mathbf{v}$, and $\rho^2 + \sigma^2 =1$.

   Let $T(Q)$ be the hyperboloid  
   centered at $\mathbf{c}$, with eccentricity $\epsilon$ and  parameters $c$ and $a$, but with  the
   axis vector   $ T(\mathbf{v}) = \mathbf{v}' $.
   That is, $T(Q)$ is the rotation about $\mathbf{c}$ of $Q$ from the axis vector $\mathbf{v}$ to the axis vector $\mathbf{v}'$, and $T(Q) = \{ \mathbf{x}: (\mathbf{x} - \mathbf{c})^T[I -  \epsilon^2 \mathbf{v}' \mathbf{v}'^T ](\mathbf{x} - \mathbf{c}) = a^2-c^2 \}$.
   Then $Q$ and $T(Q)$ are identical except for their orientation along the axis vectors $\mathbf{v}$ and $\mathbf{v}'$ respectively.
    Substituting $\mathbf{v}' =  \rho\boldsymbol{\varepsilon}_1 + \sigma \boldsymbol{\varepsilon}_n$  into the expression for $T(Q)$ and expanding yields:
    \begin{equation} \begin{split}
    T(Q) =   \{ \mathbf{x}:  (1 - \epsilon^2\rho^2)(x_1-c_1)^2 - 2\epsilon^2\rho \sigma(x_1-c_1)(x_{n}-c_n) + \sum_{j=2}^{n-1}(x_j-c_j)^2  \\
   +  (1 - \epsilon^2\sigma^2)(x_n -c_n)^2  = a^2 - c^2  \}.
  \end{split} \label{A6}  \end{equation} 
  The hyperplane $T(H\!P)$ is the rotation about $\mathbf{c}$ of $H\!P$ from the normal vector $\mathbf{h}$ to the normal vector $\boldsymbol {\varepsilon}_n$.
  Then $T(H\!P)$ passes through the point   $\mathbf{c} + \hat{h}\boldsymbol {\varepsilon}_n$, and  $T(H\!P) = \{\mathbf{x}: \boldsymbol {\varepsilon}_n(\mathbf{x} - \mathbf{c})  =\hat{h}  \}$.
  The hyperplanes $H\!P$ and $T(H\!P)$ are identical except for their orientation corresponding to the normal vectors $\mathbf{h}$ and $\boldsymbol{\varepsilon}_n$ respectively.    
  
   The hyperboloid $Q$ is related to the hyperplane $H\!P$ in the same way as the hyperboloid $T(Q)$ is related to the hyperplane $T(H\!P)$. 
   In particular, the intersection $Q \cap H\!P$ is identical to the intersection $T(Q) \cap T(H\!P)$, except for having orientations along different vectors.
    The intersection  $T(Q) \cap T(H\!P)$ is obtained by substituting $x_{n} - c_n = \hat{h}$ into the expression \eqref{A6} which yields 
     the following expression in the variables $(x_1,\ldots, x_{n-1})$ with axes parallel to the coordinate axes:
   \begin{equation}
   T(Q) \cap T(H\!P) = \{ (x_1, \ldots, x_{n-1}, c_n+\hat{h}):  \frac{(x_1 - \hat{c}_1)^2}{\hat{a}^2} + \frac{\sum_{j=2}^{n-1}(x_j - \hat{c}_j)^2}{\tilde{b}}  =1 \}, \label{tqithp}
   \end{equation}
where 
 \begin{equation} \begin{split}
 &\hat{c}_1 = c_1 + \tilde{c}, \;\;  \tilde{c} = \frac{\epsilon^2 \rho \sigma \hat{h}}{1 - \epsilon^2 \rho^2}, \;\; \hat{c}_j =c_j \;\; \text{for}\; \;j=2,\ldots,n-1   \\
 &\hat{a}^2  =\frac{(1 - \epsilon^2)[a^2(1-\epsilon^2\rho^2) - \hat{h}^2]}{(1 - \epsilon^2\rho^2)^2} \\
 &\tilde{b} = \hat{a}^2(1 - \epsilon^2 \rho^2).
 \end{split} \label{intpar} \end{equation}
 
  Let $C$ be the cone that is the asymptotic approximation of $H$, and consider the projection of $C$ and $H$ onto the affine plane 
 determined by $\mathbf{v}$ and $\mathbf{h}$ through the point $\mathbf{c}$, denoted by aff$(\mathbf{v}, \mathbf{h}, \mathbf{c})$.  
 Figure 1 illustrates the vectors $\mathbf{v}$ and $\mathbf{g}_1$ in aff$(\mathbf{v}, \mathbf{h}, \mathbf{c})$, and
the asymptotes of $C$ represented as dashed lines through the point $\mathbf{c}$.
The angle between a projected asymptote and $\mathbf{v}$ is  $\alpha$ (or $-\alpha$) where $\cos(\alpha) = a/c$.  

The projection of the hyperplane $H\!P$ onto aff$(\mathbf{v}, \mathbf{h}, \mathbf{c})$, is some line parallel to $\mathbf{g}_1$.
The the angle between $H\!P$ and $\mathbf{v}$ equals the angle between $\mathbf{g}_1$ and $\mathbf{v}$,  
which is $\beta$ (or $-\beta$) where $\rho = \mathbf{g}_1\mathbf{v} =\cos(\beta)$.

Thus, $\epsilon \rho > 1$ if and only if $\rho > a/c$ if and only if $\cos(\beta) > \cos(\alpha)$ if and only if $\beta > \alpha$, that is, if and only if the 
hyperplane $H\!P$ intersects the cone $C$, and hence the hyperboloid $H$.
Similarly,  $\epsilon \rho < 1$ if and only if the  hyperplane $H\!P$ does not intersect the cone $C$, and hence does not intersect the  hyperboloid $H$.
Also, $\epsilon \rho = 1$ if and only if the  hyperplane $H\!P$ is parallel to an asymptotic ray of the cone $C$, so that $H\!P$ intersects either sheet of $C$ 
or is coincident with the surface of $C$ on both sheets. 

Suppose $\epsilon \rho > 1$.  Then  $HP \cap Q \ne \emptyset$.
Equation (\ref{intpar}), with $\epsilon > 1$, implies  $\hat{a}^2 > 0$. 
Also,  $\tilde{b} < 0$. 
Let $\hat{b}^2 = -\tilde{b}$.
Then $T(Q) \cap T(H\!P)$ is written as  $\frac{(x_1 - \hat{c}_1)^2}{\hat{a}^2} - \frac{\sum_{j=2}^{n-1}(x_j - \hat{c}_j)^2}{\hat{b}^2}  =1 $
which is a hyperboloid with axes parallel to the coordinate axes  $\boldsymbol{\varepsilon}_j$, for $j = 1,\ldots,n-1$.
 Therefore,  $Q \cap H\!P$ is a hyperboloid of dimension $n-1$ if and only if $T(Q) \cap T(H\!P)$ is a hyperboloid of dimension $n-1$, if and only if $\epsilon \rho > 1$ and $\hat{h}^2 \geq a^2(1-\rho^2\epsilon^2)$.  
 The eccentricity of $T(Q) \cap T(H\!P)$ is $\sqrt{\hat{a}^2 + \hat{b}^2}/\hat{a} = \sqrt{\hat{a}^2 -\hat{a}^2(1 - \epsilon^2 \rho^2)}/\hat{a} = \epsilon \rho$.

  Suppose $\epsilon \rho < 1$, and $\hat{h}^2 \geq a^2(1-\rho^2\epsilon^2))$.  Then $\epsilon > 1$  implies that $\hat{a}^2 > 0$, and $\tilde{b} > 0$. Let $\hat{b}^2 = \tilde{b}$.
  Then $T(Q) \cap T(H\!P)$ is written as  $\frac{(x_1 - \hat{c}_1)^2}{\hat{a}^2} + \frac{\sum_{j=2}^{n-1}(x_j - \hat{c}_j)^2}{\hat{b}^2}  =1 $
which is an ellipsoid with axes parallel to the coordinate axes  $\boldsymbol{\varepsilon}_j$, for $j = 1,\ldots,n-1$.
  The value $a\sqrt{1-\rho^2\epsilon^2}$ is the minimum distance between the hyperplane $H\!P$ through the center $\mathbf{c}$ 
  and any point on either sheet of the hyperboloid $Q$, so that $\hat{h}^2 \geq a^2(1-\rho^2 \epsilon^2)$ guarantees that $H\!P \cap Q \neq \emptyset$.  
Therefore,  $Q \cap H\!P$ is an ellipsoid of dimension $n-1$ if and only if $T(Q) \cap T(H\!P)$ is an ellipsoid of dimension $n-1$, if and only if  
$\epsilon \rho < 1$ and $\hat{h}^2 \geq a^2(1-\rho^2\epsilon^2))$. 
The eccentricity of $T(Q) \cap T(H\!P)$ is $\sqrt{\hat{a}^2 - \hat{b}^2}/\hat{a} = \sqrt{\hat{a}^2 -\hat{a}^2(1 - \epsilon^2 \rho^2)}/\hat{a} = \epsilon \rho$.  
  
Suppose $\epsilon \rho = 1$, so that $Q \cap H\!P \neq \emptyset$.
With $x_n - c_n = \hat{h}$, expression \eqref{A6} becomes  
\begin{equation}
T(Q) \cap T(H\!P) = \{ (x_1,\ldots,x_{n-1}, c_n+\hat{h}):  4\tilde{c}(x_1 - \hat{c}_1) = \sum_{j=2}^{n-1}(x_j-\hat{c}_j)^2  \} \label{tqithppar}
\end{equation}
where 
\begin{equation}
\tilde{c} = \frac{\epsilon \sigma \hat{h}}{2},\;\;\;\; \hat{c}_1 = c_1 + \frac{(1- \epsilon^2 \sigma^2)\hat{h}^2 - (a^2 - c^2)}{4\tilde{c}}   ,\;\;  \hat{c}_j = c_j \;\; \text{for}\;\; j = 2,\ldots,n-1  \label{para}
\end{equation}
which is a paraboloid of dimension $n-1$ with axes parallel to the coordinate axes.  Thus 
 $Q \cap H\!P$ is a paraboloid of dimension $n-1$ if and only if  $T(Q) \cap T(H\!P)$ is a paraboloid of dimension $n-1$, if  $\epsilon \rho = 1$.   \hfill $\Box$
 
The following corollary gives the expressions to compute the vectors and parameters of the $n-1$ dimensional conic section resulting from the intersection of a hyperplane and  
an $n$-dimensional conic section.

\begin{corollary} If $H\!P \cap Q$ is a hyperboloid or an ellipsoid, then its vectors are given as follows: axis vector $\mathbf{g}_1$, center  
$\hat{\mathbf{c}} = \mathbf{c} + \hat{h}\mathbf{h} + \tilde{c}\mathbf{g}_1$, vertices $\hat{\mathbf{a}} = \hat{\mathbf{c}} \pm \hat{a}\mathbf{g}_1$, 
focal points $\hat{\mathbf{c}} \pm \hat{a}\epsilon \rho \mathbf{g}_1$, and its parameters $\hat{a}$ and $\hat{b}$ are given by \eqref{intpar}.  
If $H\!P \cap Q$ is a paraboloid, then its vectors are given as follows: axis vector $\mathbf{g}_1$, 
center  $\hat{\mathbf{c}} = \mathbf{c} +\hat{h}\mathbf{h} + \hat{c}\mathbf{g}_1$, and focal points $\hat{\mathbf{c}} \pm \tilde{c} \mathbf{g}_1$, 
with  $\tilde{c}$  and $\hat{c}$ given by  \eqref{para}.
\end{corollary}

The next Corollary shows that the paths traced out by the centers and vertices of the conic section resulting from $Q \cap H\!P$ are continuous with respect to $\rho$.

\begin{corollary} For hyperboloids resulting from the intersections  $Q \cap H\!P$ with $1/\epsilon < \rho < \infty$, the paths of the centers $\hat{\mathbf{c}}$ and the vertices $\hat{\mathbf{a}}$ are continuous.    For ellipsoids resulting from the intersections  $Q \cap H\!P$ with $0 <  \rho < 1/\epsilon$, the paths  of the centers $\hat{\mathbf{c}}$ and the vertices $\hat{\mathbf{a}}$ are continuous.  The path of the centers $\hat{\mathbf{c}}$ and the vertices $\hat{\mathbf{a}}$  are continusous at $\rho = 1/\epsilon$ when $Q \cap H\!P$ is a paraboloid.
\end{corollary}
\noindent \textbf{Proof:} 
The expressions for $\tilde{c}$, $\hat{a}$ and $\hat{b}$, in \eqref{intpar} and Corollary 12, are used to show that the centers $\hat{\mathbf{c}}$ and vertices $\hat{\mathbf{a}}$ of $Q \cap H\!P$ are continuous for $1/\epsilon < \rho < \infty$, and $0 < \rho < 1/\epsilon$.  
To show that $\lim_{\rho \rightarrow 1/\epsilon}$ of \eqref{tqithp} yields \eqref{tqithppar}, first observe that 
 $\lim_{\rho \rightarrow 1/\epsilon} \frac{(\epsilon^2\rho \sigma \hat{h})^2}{(1- \epsilon^2\rho^2)} = \epsilon^2 \sigma^2\hat{h}^2$.
Expression \eqref{tqithp} is re-written as
\begin{equation}
\frac{(1- \epsilon^2\rho^2)(x_1-c_1)^2 + 2\epsilon^2\rho \sigma \hat{h}(x_1 - c_1) + \frac{(\epsilon^2\rho \sigma \hat{h})^2}{(\epsilon^2\rho^2-1)}}{(a^2 - c^2) -(1- \epsilon^2 \sigma^2 )\hat{h}^2 +  \frac{(\epsilon^2\rho \sigma \hat{h})^2}{(1- \epsilon^2\rho^2)}} - \frac{\sum_{j=2}^{n-1}(x_j - c_j)^2}{(a^2 - c^2) -(1- \epsilon^2 \sigma^2 )\hat{h}^2 +  \frac{(\epsilon^2\rho \sigma \hat{h})^2}{(1- \epsilon^2\rho^2)}} = 1 \notag
\end{equation}
so that the limit of the left hand side  as $\rho \rightarrow 1/\epsilon$ simplifies to expression \eqref{tqithppar}. \hfill $\Box$

The intersection of a hyperplane and a cone is a special case of the intersection of a hyperplane with  hyperboloid, and is characterized by the following corollary.
\begin{corollary}  Suppose $C = \{ \mathbf{x}: (\mathbf{x} - \mathbf{c})^T[I -  \epsilon^2 \mathbf{v} \mathbf{v}^T ] (\mathbf{x} - \mathbf{c}) = 0  \}$ 
 is a cone in $\real^n$,  centered at $\mathbf{c} = (c_1, \ldots, c_n)^T$, 
with  axis  vector $\mathbf{v}$ of unit length,  eccentricity $\epsilon$, and 
  parameters $a$ and $c$, and suppose 
  $H\!P = \{ \mathbf{x}: \mathbf{h}( \mathbf{x} -\mathbf{c})=\hat{h} \}$,
  is a hyperplane with  $\| \mathbf{h} \| =1$ and $\hat{h} > 0$.
  Let $\rho = \sqrt{1-(\mathbf{h} \mathbf{v})^2}$. 
   Then $C \cap H\!P$  is a hyperboloid (ellipsoid, paraboloid) of dimension $n-1$ iff  $\epsilon \rho > 1$ ($\epsilon \rho < 1$, $\epsilon \rho = 1$). 
 \end{corollary}
  
  The parameters and vectors of $C \cap H\!P$ are determined by the same expressions as for the intersection  $Q \cap H\!P$, but without the term $(c^2 -a^2)$ in  \eqref{intpar} and \eqref{para}.
  Since $\hat{h} > 0$, the hyperplane $H\!P$ will always intersect either one sheet or both sheets of the cone.
    
    \begin{figure}
\begin{center}
\includegraphics[scale=.7]{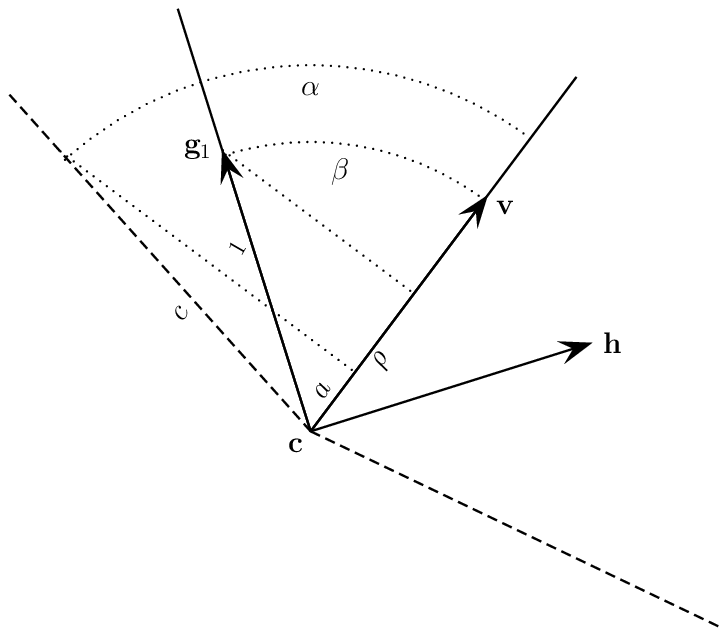}
\caption{Hyperbola $H$ projected onto the affine plane aff$(\mathbf{v}, \mathbf{h}, \mathbf{c})$}
\end{center}
\end{figure}

  Then $\epsilon \rho  >(<,=) \;1$ if and only if $\rho > (<,=) \;1/\epsilon = a/c$ if and only if $\cos( \beta) >(<,=) \cos(\alpha)$ if and only $\beta <(>,=)\; \alpha$, that is, if and only the vector $\mathbf{g}_1$ is inside (outside, parallel) to the asymptotes of the cone.


\section{Intersections of a hyperplane with an ellipsoid}

It is well known that the nonempty intersection of a hyperplane with an ellipsoid is always an ellipsoid of dimension $n-1$.  
This result is stated below using the representation $Q$ for an ellipsoid. 
The resulting parameters and vectors are analogous to the intersection of a hyperplane with a hyperboloid, but with $c < a$ in  representation \eqref{HQ1}.

\begin{property}  Suppose $Q = \{ \mathbf{x}: (\mathbf{x} - \mathbf{c})^T[I -  \epsilon^2 \mathbf{v} \mathbf{v}^T ] (\mathbf{x} - \mathbf{c}) = a^2 -c^2  \}$ 
 is an ellipsoid in $\real^n$,  centered at $\mathbf{c} = (c_1, \dots, c_n)^T$, with  axis  vector $\mathbf{v}$ of unit length,  eccentricity $\epsilon < 1$, 
 so that  $a > c$, and suppose $H\!P = \{ \mathbf{x}: \mathbf{h}( \mathbf{x} -\mathbf{c})=\hat{h} \}$,
  is a hyperplane  with  $\| \mathbf{h} \| =1$.
  Let $\rho = \sqrt{1-(\mathbf{h} \mathbf{v})^2}$.  
 Then $Q \cap H\!P$  is an ellipsoid of dimension $n-1$ if $\hat{h}^2 < a^2(1- \epsilon^2\rho^2)$.
 \end{property}
 \noindent \textbf{Proof:}  Applying the transformation $T$ used in the proof of Property 6.1, to $Q$ and to $H\!P$  leads to expressions  \eqref{tqithp} and \eqref{intpar} for $T(Q) \cap T(H\!P)$, but with $c < a$ and $\epsilon^2\rho^2 <1$, since $\epsilon<1$ and $\rho \leq 1$. 
 
 The value of $a\sqrt{1- \epsilon^2\rho^2}$ is the maximum distance between a point on the ellipsoid $Q$ and the hyperplane $H\!P$ through the center $\mathbf{c}$, so that $\hat{h}^2 < a^2(1- \epsilon^2\rho^2)$ guarantees that $H\!P \cap Q \neq \emptyset$.  
 Since $\hat{h}^2 < a^2(1- \epsilon^2\rho^2)$,  and $\epsilon < 1$, $\hat{a}^2 > 0$, and 
  $\tilde{b}> 0$. Let $\hat{b}^2 = \tilde{b}$. 
   Then $T(Q) \cap T(H\!P)$ is written as  $\frac{(x_1 - \hat{c}_1)^2}{\hat{a}^2} + \frac{\sum_{j=2}^{n-1}(x_j - \hat{c}_j)^2}{\hat{b}^2}  =1 $
which is an ellipsoid with axes parallel to the coordinate axes  $\boldsymbol{\varepsilon}_j$, for $j = 1,\ldots,n-1$.
Therefore,  $Q \cap H\!P$ is an ellipsoid of dimension $n-1$ if and only if $T(Q) \cap T(H\!P)$ is an ellipsoid of dimension $n-1$, if  $\epsilon \rho < 1$ and $\hat{h}^2 \geq a^2(1-\rho^2\epsilon^2))$. 
The eccentricity of $T(Q) \cap T(H\!P)$ is $\sqrt{\hat{a}^2 - \hat{b}^2}/\hat{a} = \sqrt{\hat{a}^2 -\hat{a}^2(1 - \epsilon^2 \rho^2)}/\hat{a} = \epsilon \rho$.  
The vectors of $Q \cap H\!P$ are: center  
$\hat{\mathbf{c}} = \mathbf{c} + \hat{h}\mathbf{h} -\frac{\epsilon^2 \rho \sigma \hat{h}}{1 - \epsilon^2 \rho^2}\mathbf{g}_1$, vertices $\hat{\mathbf{a}} = \hat{\mathbf{c}} \pm \hat{a}\mathbf{g}_1$, and focal points $\hat{\mathbf{c}} \pm \hat{a}\epsilon \rho \mathbf{g}_1$. \hfill $\Box$


\section{Intersections of a hyperplane with a paraboloid}

The intersection of a hyperplane with a paraboloid results in either a paraboloid of dimension $n-1$, or an ellipsoid of dimension $n-1$.  

\begin{property} Suppose $P = \{  \mathbf{x}:  ( \mathbf{x} - \mathbf{c})^T [I -  \mathbf{v} \mathbf{v}^T] ( \mathbf{x} - \mathbf{c})  = 4c\mathbf{v}(\mathbf{x} - \mathbf{c}) \}$, is a paraboloid in $\real^n$, centered at $\mathbf{c} = (c_1,\ldots,c_n)^T$, with axis vector $\mathbf{v}$ of unit length, and parameter $c$, and suppose that   $H\!P = \{ \mathbf{x}: \mathbf{h}( \mathbf{x} -\mathbf{c}) =\hat{h} \}$,
  is a hyperplane where  $\| \mathbf{h} \| =1$.
  Suppose that  $P \cap H\!P$ is nonempty.
  Let $\rho = \sqrt{1-(\mathbf{h} \mathbf{v})^2}$.  
 Then $P \cap H\!P$  is a paraboloid of dimension $n-1$ if $\rho = 1$, or an ellipsoid of dimension $n-1$ if $\rho < 1$.
 \end{property}
 \noindent \textbf{Proof:} 
 Assume $\rho = 1$. 
 Applying the transformation $T$ used in the proof of Property 6.1, to $P$ and to $H\!P$, gives $\mathbf{h}\mathbf{v} = \sigma = 0$ so that   $\mathbf{g}_1 = \mathbf{v}$, and $T(\mathbf{v}) = \mathbf{v}' = \boldsymbol{\varepsilon}_1$.
  Let $T(P)$ be the paraboboloid  
   centered at $\mathbf{c}$, with parameter  $c$, but with  the
   axis vector   $ T(\mathbf{v}) = \mathbf{v}' $.
   Then $P$ and $T(P)$ are identical except for their orientation along the axis vectors $\mathbf{v}$ and $\mathbf{v}'$ respectively.
    Substituting $\mathbf{v}' =  \boldsymbol{\varepsilon}_1$  into the expression for $P$ and expanding yields:
    \begin{equation} \begin{split}
    T(P) =   \{ \mathbf{x}:  \sum_{j=2}^{n}(x_j-c_j)^2  = 4c(x_1 - c_1) \}.
  \end{split} \label{TP}  \end{equation} 
 
   The hyperplane  $T(H\!P) = \{\mathbf{x}: \boldsymbol {\varepsilon}_n(\mathbf{x} - \mathbf{c})  =\hat{h}  \}$
  passes through the point $\mathbf{c} + \hat{h}\boldsymbol{\varepsilon}_n$, and has normal vector $\boldsymbol{\varepsilon}_n$.
  Then $H\!P$ and $T(H\!P)$ are identical except for their orientation corresponding to the normal vectors $\mathbf{h}$ and $\boldsymbol{\varepsilon}_n$ respectively.    
   The paraboloid $P$ is related to the hyperplane $H\!P$ in the same way as the paraboloid $T(P)$ is related to the hyperplane $T(H\!P)$. 
   In particular, the intersection $P \cap H\!P$ is equivalent to the intersection $T(P) \cap T(H\!P)$.
    The intersection  $T(P) \cap T(H\!P)$ is obtained by substituting $x_{n} - c_n = \hat{h}$ into the expression \eqref{TP} which yields 
     the following expression in the variables $(x_1,\ldots, x_{n-1})$:
    $$T(Q) \cap T(H\!P) = \{ (x_1, \ldots, x_{n-1}, c_n+\hat{h}) :   \sum_{j=2}^{n-1}(x_j - c_j)^2   = 4c(x_1 - \hat{c}_1)  \},$$
    where $\hat{c}_1 = c_1 +\tilde{c}$ and $\tilde{c} = \hat{h}^2 / 4c$.
    Thus $T(Q) \cap T(H\!P)$ is a paraboloid with axes parallel to the coordinate axes, and $P \cap H\!P$ is a paraboloid of dimension $n-1$.  
    
    If $\rho < 1$, then $\sigma > 0$  and $T(\mathbf{v}) = \mathbf{v}' = \rho \boldsymbol{\varepsilon}_1 + \sigma \boldsymbol{\varepsilon}_n $.  Thus 
     \begin{equation} \begin{split}
    T(P) =   \{ \mathbf{x}:  (\rho^2 -1)(x_1-c_1)^2 + 2\rho \sigma(x_1-c_1)(x_{n}-c_n) - \sum_{j=2}^{n-1}(x_j-c_j)^2 - \\
   +  (\sigma^2 -1)(x_n -c_n)^2  = -4c \rho(x_1 - c_1) -4c \sigma (x_n-c_n)  \}.
  \end{split} \label{parapara}  \end{equation} 
  The intersection  $T(P) \cap T(H\!P)$ is obtained by substituting $x_{n} - c_n = \hat{h}$ into the expression \eqref{TP} which yields 
     the following expression in the variables $(x_1,\ldots, x_{n-1})$:

$$T(P) \cap T(H\!P) = \{ (x_1, \ldots, x_{n-1}, c_n+\hat{h}) : \frac{(x_1-\hat{c}_1)^2}{\hat{a}^2} - \frac{\sum_{j=2}^{n-1}(x_j - \hat{c}_j)^2}{\hat{b}^2}   = 1 \},$$
where
\begin{equation}\begin{split}
\hat{c}_1 = c_1- \tilde{c}, \;\; \tilde{c} = \frac{\rho(\sigma\hat{h} + 2c)}{\rho^2-1}, \;\; \hat{c}_j = c_j \;\;\text{for}\;\;j-2,\ldots,n-1, \\
 \hat{a}^2 = \frac{4c(c\rho^2+\sigma\hat{h})}{(\rho^2-1)^2}, \;\;\;\text{and} \;\;\; \hat{b}^2 = \hat{a}^2(\rho^2-1).
\end{split} \label{raraellip} \end{equation}
With $\rho < 1$, and if $\hat{h} > -c\rho^2/\sigma$, the coefficient $\hat{a}^2 >0$.  The expression $|-c\rho^2/\sigma|$ is the maximum distance between the paraboloid and the hyperplane $H\!P$ through the center $\mathbf{c}$, so that $\hat{h} > -c\rho^2/\sigma$ guarantees that $H\!P \cap P \neq \emptyset$.  With $\hat{a}^2 >0$, the  coefficient $\frac{-1}{\hat{b}^2} >0$ for each of the remaining terms, so that $T(P) \cap T(H\!P)$ is an ellipsoid of dimension $n-1$ with axes parallel to the coordinate axes.  Thus $P \cap H\!P$ is an ellipsoid of dimension $n-1$.\hfill $\Box$

   For the case where $\rho = 1$, the vectors and parameters of the paraboloid $P \cap H\!P$ are given as follows: the axis vector is $\mathbf{g}_1$, and 
    the center (and vertex) is given by  $\hat{\mathbf{c}} = \mathbf{c} + \hat{h}\mathbf{h} + \frac{\hat{h}^2}{4c} \mathbf{g}_1$.
    For the case where $\rho < 1$, the vectors and parameters of the ellipsoid $P \cap H\!P$ are given as follows: the axis vector is $\mathbf{g}_1$, 
    the center is given by  $\hat{\mathbf{c}} = \mathbf{c}  + \hat{h}\mathbf{h} - \tilde{c}\mathbf{g}_1$,  the vertices are given by $\hat{\mathbf{c}}   \pm \hat{a} \mathbf{g}_1$, and the focal points by $\hat{\mathbf{c}}   \pm \hat{a}\epsilon \rho \mathbf{g}_1$, 
 
 
 \section{A class of hyperboloids}
  
  A class of hyperboloids is defined whose pair-wise intersections have  additional properties.
Let $\mathcal{P} = \{ \mathbf{p}_1, \ldots, \mathbf{p}_m \}$ be a finite set of  points in $\real^n$, and for each  point $\mathbf{p}_i \in \mathcal{P}$ 
let $r_i$ be a non-negative radius and let $[\mathbf{p}_i, r_i] = \{ \mathbf{x}: \parallel \mathbf{p}_i - \mathbf{x} \parallel \leq r_i \}$ denote the corresponding Euclidean ball.
For each pair of balls $[\mathbf{p}_j, r_j]$ and $[\mathbf{p}_k, r_k]$, the \textbf{bisector} $B_{j,k}$ of $[\mathbf{p}_j, r_j]$ and $[\mathbf{p}_k, r_k]$ is the set of points 
$\mathbf{x}$ whose distance to $\mathbf{p}_j$ plus $r_j$ equals the distance to $\mathbf{p}_k$ plus $r_k$.  That is
\begin{equation}
B_{j,k} = \{ \mathbf{x}: \parallel \mathbf{p}_j - \mathbf{x} \parallel + r_j = \parallel \mathbf{p}_k - \mathbf{x} \parallel + r_k \}. \label{bisec} \end{equation}
For each point $\mathbf{x} \in B_{j,k}$, the ball $[\mathbf{x}, z]$, with radius $z = \parallel \mathbf{p}_j - \mathbf{x} \parallel + r_j = \parallel \mathbf{p}_k - \mathbf{x} \parallel + r_k$, contains the balls  $[\mathbf{p}_j, r_j]$ and $[\mathbf{p}_k, r_k]$ and is internally tangent to each.

If $r_j = r_k$, the bisector $B_{j,k}$ is the hyperplane that bisects, and is orthogonal to, the line segment between $\mathbf{p}_j$ and $\mathbf{p}_k$, that is,
\begin{equation}
B_{j,k} = \{ \mathbf{x}: (\mathbf{p}_j - \mathbf{p}_k)\mathbf{x} = .5(\mathbf{p}_j - \mathbf{p}_k)(\mathbf{p}_j + \mathbf{p}_k) \}. \label{hyp} \end{equation}

However, if $r_j \neq r_k$, and choosing $j$ and $k$ so that $r_j > r_k$, expression \eqref{H} defines a hyperboloid wth focal points $\mathbf{p}_j$ and  $\mathbf{p}_k$  and constant  $2a_{j,k} = r_j - r_k$, provided   $a_{j,k} < c_{j,k}$.   The bisector $B_{j,k} = H_j$ is the sheet, defined by \eqref{H1},  
that is, $B_{j,k} = \{ \mathbf{x}: \parallel \mathbf{p}_k - \mathbf{x} \parallel  - \parallel \mathbf{p}_j - \mathbf{x} \parallel = r_j - r_k \}.$
The class of hyperboloids is determined by all pairs of points $\mathbf{p}_j, \mathbf{p}_k \in \mathcal{P}$ with $r_j > r_k$ and $a_{j,k} < c_{j,k}$.

The condition $a_{j,k} < c_{j,k}$, is equivalent to $r_j - r_k = 2a_{j,k} < \parallel \mathbf{p}_j - \mathbf{p}_k \parallel $, which is equivalent to the condition that  neither one of the two balls $[\mathbf{p}_j,r_j]$ and $[\mathbf{p}_k, r_k]$ is contained in the other.

The vertex $\mathbf{a}_{j,k}$ of the sheet $B_{j,k}$ is the center of the ball with minimum radius that contains the balls $[\mathbf{p}_j,r_j]$ and $[\mathbf{p}_k,r_k]$.  If $r_j = r_k$, the point  $\mathbf{c}_{j,k} = (\mathbf{p}_j + \mathbf{p}_k)/2$ is the center of the ball with minimum radius that contains the balls  $[\mathbf{p}_j,r_j]$ and $[\mathbf{p}_k,r_k]$.  
Each point $\mathbf{x}$ on the opposite sheet is the center of a ball that is externally tangent to the two balls  $[\mathbf{p}_j,r_j]$ and $[\mathbf{p}_k,r_k]$.

Let  $T = \{ \mathbf{p}_{j}, \mathbf{p}_{k}, \mathbf{p}_{l} \}$ be a set of three affinely independent
points from $\mathcal{P}$ ordered so that $r_{j} \geq r_{k} \geq r_{l}$, and let $\mathcal{B} = \{ B_{j,k}, B_{j,l} , B_{k,l} \}$
denote the bisectors corresponding to the respective pairs of points from $T$.
Each bisector in $\mathcal{B}$ may be either a sheet of a hyperboloid or a hyperplane.
If the radii are unequal, that is  $r_{j} > r_{l} $,
then at least two pairs of points from $T$ have unequal radii, 
and the bisector corresponding to each of these pairs is a sheet of the corresponding hyperboloid. 
If $r_j = r_l$, each of the three bisectors is  a hyperplane.

The following theorem constructs a unique hyperplane $H_T$  that contains the intersection of  any two 
bisectors in $\cal{B}$.  The theorem also shows that for any two bisectors in $\cal{B}$, say $B_{j,k}$ and $B_{j,l}$, that
$ B_{j,k} \cap B_{j,l} =  B_{j,k}\cap H_T = B_{j,l}\cap H_T$.
Also, if one of the bisectors in $\cal{B}$ is a hyperplane, then it is identical to $H_T$.
This result allows the intersection of two bisectors to be determined by the intersection of either one of the bisectors with the hyperplane $H_T$.
This result leads to a procedure for finding the intersection of bisectors determined by all pairs of points in a subset of  $\mathcal{P}$.

Theorem 17 extends a result in reference \cite{leva}  that assumes only  two hyperboloids with a common focal point.
 
\begin{theorem}
Suppose that   $T = \{ \mathbf{p}_{j}, \mathbf{p}_{k}, \mathbf{p}_{l} \}$ is a subset of three affinely independent points from $\mathcal{P}$,
   ordered so that $r_{j} \geq r_{k} \geq r_{l}$, with $r_{j} > r_{l}$,
   and suppose that the intersection of the three bisectors from $\mathcal{B} = \{ B_{j,k}, B_{j,l} , B_{k,l} \}$ is nonempty.  
   Then there exist the hyperplane $H_T = \{ \mathbf{x}: \mathbf{h}_T \mathbf{x}  =  d_T \}$ so that for each pair of bisectors in $\cal{B}$, their intersection equals the intersection of $H_T$ with either bisector in the pair, and if any one of the bisectors in $\cal{B}$ is a hyperplane, it is identical to $H_T$.    
   \end{theorem}
   \noindent \textbf{Proof:}:
  Suppose first that the radii corresponding to the points in $T$ satisfy   $r_j > r_k > r_l$, so that each of the three bisectors in $\cal{B}$ is a sheet of a hyperboloid.
  First, consider the pair of bisectors $B_{j,k}$ and $B_{j,l}$ from $\cal{B}$.
  Since the vectors $\mathbf{v}_{j,k}$ and $\mathbf{v}_{j,l}$ are linearly independent,  there exists a point $\mathbf{d}_{jkjl}$
  that  is the intersection point of the directrix $\mathbf{v}_{j,k} \mathbf{x} = \mathbf{v}_{j,k} \mathbf{d}_{j,k}$ of $B_{j,k}$, the directrix
  $\mathbf{v}_{j,l} \mathbf{x} = \mathbf{v}_{j,l} \mathbf{d}_{j,l}$ of $B_{j,l}$, and  aff$(T)$.
  Thus $\mathbf{v}_{j,k} \mathbf{d}_{j,k} = \mathbf{v}_{j,k} \mathbf{d}_{jkjl}$ and $\mathbf{v}_{j,l} \mathbf{d}_{j,l} = \mathbf{v}_{j,l} \mathbf{d}_{jkjl}$. 
   Combining this result with equation \eqref{HH1} yields  
$B_{j,k} = \{  \mathbf{x}:   \parallel \mathbf{p}_{j} - \mathbf{x} \parallel = \epsilon_{j,k} \mathbf{v}_{j,k} \mathbf{x} -  \epsilon_{j,k}  \mathbf{v}_{j,k} \mathbf{d}_{jkjl} \} $ and
$B_{j,l} = \{  \mathbf{x}:   \parallel \mathbf{p}_{j} - \mathbf{x} \parallel = \epsilon_{j,l} \mathbf{v}_{j,l} \mathbf{x} -  \epsilon_{j,l}   \mathbf{v}_{j,l} \mathbf{d}_{jkjl}  \} $. 
If $\mathbf{x} \in B_{j,k} \cap B_{j,l}$, then $\mathbf{x}$  satisfies
$(\epsilon_{j,k} \mathbf{v}_{j,k} - \epsilon_{j,l} \mathbf{v}_{j,l}) \mathbf{x} =  (\epsilon_{j,k}  \mathbf{v}_{j,k} -  \epsilon_{j,l}   \mathbf{v}_{j,l}) \mathbf{d}_{jkjl}$.
The hyperplane $H_T = \{ \mathbf{x}: \mathbf{h}_T\mathbf{x} = d_T \}$ is defined as:
\begin{equation}
H_T = \{ \mathbf{x}: (\epsilon_{j,k} \mathbf{v}_{j,k} - \epsilon_{j,l} \mathbf{v}_{j,l}) \mathbf{x} =  (\epsilon_{j,k}  \mathbf{v}_{j,k} -  \epsilon_{j,l}   \mathbf{v}_{j,l}) \mathbf{d}_{jkjl} \} \label{HT}
\end{equation}
where $\mathbf{h}_T = (\epsilon_{j,k} \mathbf{v}_{j,k} - \epsilon_{j,l} \mathbf{v}_{j,l})$ 
is the normal vector, and $d_T = \mathbf{h}_T \mathbf{d}_{jkjl}$
is the right hand side value.
Then $B_{j,k} \cap B_{j,l} \subseteq H_T \cap B_{j,k}$.
Conversely, if $\mathbf{x} \in H_T \cap B_{j,k}$, then $\mathbf{x}$ satisfies 
$\epsilon_{j,k} \mathbf{v}_{j,k} \mathbf{x} -  \epsilon_{j,k}  \mathbf{v}_{j,k} \mathbf{d}_{jkjl} = 
\epsilon_{j,l} \mathbf{v}_{j,l} \mathbf{x} -  \epsilon_{j,l}   \mathbf{v}_{j,l} \mathbf{d}_{jkjl}$
and $ \parallel \mathbf{p}_{j} - \mathbf{x} \parallel = \epsilon_{j,k} \mathbf{v}_{j,k} \mathbf{x} -  \epsilon_{j,k}  \mathbf{v}_{j,k} \mathbf{d}_{jkjl}$.
Combining these equations shows that $\mathbf{x}$ satisfies  
 $\parallel \mathbf{p}_{j} - \mathbf{x} \parallel = \epsilon_{j,l} \mathbf{v}_{j,l} \mathbf{x} -  \epsilon_{j,l}   \mathbf{v}_{j,l} \mathbf{d}_{jkjl}$,
 so that $\mathbf{x} \in B_{j,l}$. Thus $ H_T \cap B_{j,k} = B_{j,k} \cap B_{j,l}$.
 A parallel argument shows that  $ H_T \cap B_{j,l} = B_{j,k} \cap B_{j,l}$.

Next, consider the pair of bisectors  $B_{j,k}$ and $B_{k,l}$, and let $\mathbf{d}_{jkkl}$ be the point of intersection 
  of the directrix $\mathbf{v}_{j,k} \mathbf{x} = \mathbf{v}_{j,k} \mathbf{d}_{j,k}$ of $B_{j,k}$, the directrix
  $\mathbf{v}_{k,l} \mathbf{x} = \mathbf{v}_{k,l} \mathbf{d}_{k,l}$ of $B_{k,l}$, and  aff$(T)$. 
  From expression \eqref{H1} and \eqref{HH1}, the hyperboloids $B_{j,k}$ and $B_{k,l}$ may be written as 
$B_{j,k} = \{  \mathbf{x}:  \; \parallel \mathbf{p}_{j} - \mathbf{x} \parallel =  \parallel \mathbf{p}_{k} - \mathbf{x} \parallel -2a_{j,k} =
 \epsilon_{j,k} \mathbf{v}_{j,k} \mathbf{x} -  \epsilon_{j,k}  \mathbf{v}_{j,k} \mathbf{d}_{jkkl}  \} $ and
$B_{k,l} = \{  \mathbf{x}:   \parallel \mathbf{p}_{k} - \mathbf{x} \parallel = \epsilon_{k,l} \mathbf{v}_{k,l} \mathbf{x} -  \epsilon_{k,l}   \mathbf{v}_{k,l} \mathbf{d}_{jkkl}  \} $
respectively.
Thus all points in the intersection of  $B_{j,k}$ and $B_{k,l}$ must satisfy
$\epsilon_{j,k} \mathbf{v}_{j,k} \mathbf{x} -  \epsilon_{j,k}  \mathbf{v}_{j,k} \mathbf{d}_{jkkl}  +2a_{jk} = 
\epsilon_{k,l} \mathbf{v}_{k,l} \mathbf{x} -  \epsilon_{k,l}   \mathbf{v}_{k,l} \mathbf{d}_{jkkl} $
which is expressed as the hyperplane 
\begin{equation}
\{ \mathbf{x}: (\epsilon_{j,k} \mathbf{v}_{j,k}  - \epsilon_{k,l} \mathbf{v}_{k,l}) \mathbf{x}  = (\epsilon_{j,k} \mathbf{v}_{j,k}  - \epsilon_{k,l} \mathbf{v}_{k,l}) \mathbf{d}_{jkkl} -2a_{jk} \}.
\label{Ajkkl} \end{equation}

To show the equivalence of  $H_T$, given by \eqref{HT}, and \eqref{Ajkkl},   the normal vector of $H_T$ is expressed as
$$ (\epsilon_{j,k} \mathbf{v}_{j,k}  - \epsilon_{j,l} \mathbf{v}_{j,l}) = \frac{r_j(\mathbf{p}_l - \mathbf{p}_k) + r_k(\mathbf{p}_j - \mathbf{p}_l) + r_l(\mathbf{p}_k - \mathbf{p}_j)}
{(r_j-r_k)(r_j-r_l)},$$ by substituting the definition of each parameter and vector, and simplifying.  Similarly, the 
normal vector of \eqref{Ajkkl} is expressed as 
$$ (\epsilon_{j,k} \mathbf{v}_{j,k}  - \epsilon_{k,l} \mathbf{v}_{k,l}) = \frac{r_j(\mathbf{p}_l - \mathbf{p}_k) + r_k(\mathbf{p}_j - \mathbf{p}_l) + r_l(\mathbf{p}_k - \mathbf{p}_j)}
{(r_j-r_k)(r_k-r_l)},$$
so that $$(\epsilon_{j,k} \mathbf{v}_{j,k}  - \epsilon_{j,l} \mathbf{v}_{j,l}) = \frac{(r_k-r_l)}{(r_j-r_l)}(\epsilon_{j,k} \mathbf{v}_{j,k}  - \epsilon_{k,l}\mathbf{v}_{k,l}).$$
Using the same approach, the right hand side of \eqref{HT} is expressed as 
$$(\epsilon_{j,k} \mathbf{v}_{j,k}  - \epsilon_{j,l} \mathbf{v}_{j,l}) \mathbf{d}_{jkjl} =
\frac{r_j(\mathbf{p}_l^2 - \mathbf{p}_k^2) + r_k(\mathbf{p}_j^2 - \mathbf{p}_l^2) + r_l(\mathbf{p}_k^2 - \mathbf{p}_j^2)}{2(r_j - r_k)(r_j - r_l)} - \frac{r_k - r_l}{2},$$
and the right hand side of \eqref{Ajkkl} is expressed as 
$$(\epsilon_{j,k} \mathbf{v}_{j,k}  - \epsilon_{j,l} \mathbf{v}_{j,l}) \mathbf{d}_{jkkl} -2a_{j,k} =
\frac{r_j(\mathbf{p}_j^2 - \mathbf{p}_k^2) + r_k(\mathbf{p}_j^2 - \mathbf{p}_l^2) + r_l(\mathbf{p}_k^2 - \mathbf{p}_l^2)}{2(r_j - r_k)(r_k - r_l)} - \frac{r_j - r_l}{2},$$
so that 
$$(\epsilon_{j,k} \mathbf{v}_{j,k}  - \epsilon_{j,l} \mathbf{v}_{j,l}) \mathbf{d}_{jkjl} =
 \frac{(r_k-r_l)}{(r_j-r_l)}[(\epsilon_{j,k} \mathbf{v}_{j,k}  - \epsilon_{k,l} \mathbf{v}_{k,l}) \mathbf{d}_{jkkl} -2a_{j,k}] .$$
 Thus the hyperplanes  \eqref{HT} and \eqref{Ajkkl} are equivalent.  A similar argument to the above shows that 
 $ H_T \cap B_{j,k} = B_{j,k} \cap B_{k,l}$ and $ H_T \cap B_{k,l} = B_{j,k} \cap B_{k,l}$.

Next consider the pair of bisectors  $B_{j,l}$ and $B_{k,l}$ and let $ \mathbf{d}_{jlkl}$ be the point of intersection of the directrix 
$\mathbf{v}_{j,l} \mathbf{x} = \mathbf{v}_{j,l} \mathbf{d}_{j,l}$ of $B_{j,l}$, the directrix
  $\mathbf{v}_{k,l} \mathbf{x} = \mathbf{v}_{k,l} \mathbf{d}_{k,l}$ of $B_{k,l}$, and  aff$(T)$.  
Then a parallel development shows that if $\mathbf{x} \in B_{j,l} \cap B_{k,l}$, then $\mathbf{x}$ must be on the hyperplane:
\begin{equation}
\{ \mathbf{x}:  (\epsilon_{j,l} \mathbf{v}_{j,l}  - \epsilon_{k,l} \mathbf{v}_{k,l}) \mathbf{x}  = (\epsilon_{j,l} \mathbf{v}_{j,l}  - \epsilon_{k,l} \mathbf{v}_{k,l}) \mathbf{d}_{jlkl} -2a_{jk} \}. 
 \label{Ajlkl} \end{equation}
  A analysis similar to the above case shows that the hyperplanes  \eqref{HT} and \eqref{Ajlkl} are equivalent, and that $ H_T \cap B_{j,l} = B_{j,l} \cap B_{k,l}$ and
  $ H_T \cap B_{k,l} = B_{j,l} \cap B_{k,l}$.
This proves the result for the case that  $r_j > r_k > r_l$.

  Next,  suppose that two points in $T$, say $\mathbf{p}_k$ and $\mathbf{p}_l$ have equal radii, so that  $r_j > r_k = r_l$.
 The bisector $B_{k,l}$ is the hyperplane 
  $(\mathbf{p}_{k} - \mathbf{p}_{l}) \mathbf{x} = (\mathbf{p}_{k} - \mathbf{p}_{l}) (\mathbf{p}_{k} + \mathbf{p}_{l})/2$,
  as given by \eqref{hyp}.
  The following development shows that the hyperplane $B_{k,l}$ is identical to the hyperplane $H_T$.
  First observe that $a_{j,k} = a_{j,l}$ since $r_{k} = r_{l}$.  
  Then the normal  vector  of $H_T$ may be written as
  $$\epsilon_{j,k} \mathbf{v}_{j,k} - \epsilon_{j,l} \mathbf{v}_{j,l}   =  \frac{\parallel \mathbf{p}_{j} - \mathbf{p}_{k}\parallel }{2 a_{j,k}} 
  \frac{ \mathbf{p}_{j} - \mathbf{p}_{k} }{ \parallel \mathbf{p}_{j} - \mathbf{p}_{k}\parallel }  - \frac{\parallel \mathbf{p}_{j} - \mathbf{p}_{l}\parallel }{2 a_{j,l}} 
  \frac{ \mathbf{p}_{j} - \mathbf{p}_{l} }{ \parallel \mathbf{p}_{j} - \mathbf{p}_{l}\parallel } \\
  = \frac{-1}{2a_{j,k}} ( \mathbf{p}_{l} - \mathbf{p}_{k} ). $$
  The right hand side of $H_T$ becomes:\begin{alignat}{3}
(\epsilon_{j,k} \mathbf{v}_{j,k}) - \epsilon_{j,l} \mathbf{v}_{j,l}) \mathbf{d}_{jkjl} 
    &=  \epsilon_{j,k} \mathbf{v}_{j,k} \mathbf{d}_{j,k} - \epsilon_{j,l} \mathbf{v}_{j,l} \mathbf{d}_{j,l} \notag \\
 & =  \frac{c_{j,k}}{a_{j,k}}\mathbf{v}_{j,k}(\mathbf{c}_{j,k} + d_{j,k} \mathbf{v}_{j,k}) - \frac{c_{j,l}}{a_{j,l}}\mathbf{v}_{j,l} (\mathbf{c}_{j,l} + d_{j,l} \mathbf{v}_{j,l}) \notag \\
 & = \frac{\parallel \mathbf{p}_{j} - \mathbf{p}_{k}\parallel(\mathbf{p}_{j} - \mathbf{p}_{k})(\mathbf{p}_j + \mathbf{p}_k )}{2 a_{j,k}\parallel \mathbf{p}_{j} - \mathbf{p}_{k}\parallel 2} - \frac{\parallel \mathbf{p}_{j} - \mathbf{p}_{l}\parallel(\mathbf{p}_{j} - \mathbf{p}_{l})(\mathbf{p}_j + \mathbf{p}_l )}{2 a_{j,l}\parallel \mathbf{p}_{j} - \mathbf{p}_{l}\parallel 2} \notag \\
 & =\frac{1}{2a_{j,k}} \frac{( \mathbf{p}_{j}^2 - \mathbf{p}_{k}^2 - \mathbf{p}_{j}^2 + \mathbf{p}_{l}^2 )}{2} \notag \\
 & = \frac{-1}{2a_{j,k}} \frac{( \mathbf{p}_{k}^2 - \mathbf{p}_{l}^2 )}{2} \notag \\
 & =  \frac{-1}{2a_{j,k}} \frac{( \mathbf{p}_{k} - \mathbf{p}_{l} )( \mathbf{p}_{k} + \mathbf{p}_{l} )}{2}. \notag 
 \end{alignat}
Thus the hyperplane $H_T$ is  equivalent to the hyperplane $B_{k,l}$.
The former arguments are used to show that $ H_T \cap B_{j,k} = B_{j,k} \cap B_{j,l}$ and $ H_T \cap B_{j,l} = B_{j,k} \cap B_{j,l}$.

The final case assumes $r_j = r_k > r_l$, and uses analogous arguments to show that $B_{j,k}$ is equivalent to $H_T$ and that $ H_T \cap B_{j,k} = B_{j,k} \cap B_{k,l}$ and $ H_T \cap B_{k,l} = B_{j,k} \cap B_{k,l}$.
This concludes the proof.
  \hfill $\Box$
 
 The next property shows how to compute the intersection point $\mathbf{d}_{jkjl}$ under the assumptions of Theorem 9.1.
 \begin{property}
Given a set $T =  \{ \mathbf{p}_{j}, \mathbf{p}_{k}, \mathbf{p}_{l} \}$ 
of affinely independent points from $P$,
so that $r_{j} \geq r_{k} \geq r_{l}$, with $r_{j} > r_{l}$,
and the hyperplane
$H_T = \{ \mathbf{x}: \mathbf{h}_T \mathbf{x}  = \mathbf{h}_T \mathbf{d}_{jkjl} \},$
  containing the intersections of $ B_{j,k}, B_{j,l} $, and $ B_{k,l} $, then
  $$ \mathbf{d}_{jkjl} = \mathbf{d}_{jl} + \frac{\mathbf{v}_{jk}(\mathbf{d}_{jk}-\mathbf{d}_{jl})}{\mathbf{v}_{jk}\mathbf{u}_T} \mathbf{u}_T, 
\;\;\;\text{where}\;\;\; \mathbf{u}_T = \frac{\mathbf{h}_T - (\mathbf{v}_{jl} \mathbf{h}_T) \mathbf{v}_{jl}}
{\parallel \mathbf{h}_T - (\mathbf{v}_{jl} \mathbf{h}_T) \mathbf{v}_{jl}\parallel}$$ 
is a unit vector  in aff$(T)$ that is orthogonal to $\mathbf{v}_{jk}$.  
\end{property}  
\noindent \textbf{Proof:}
Direct substitution shows that 
$ \mathbf{v}_{j,k}\mathbf{d}_{jkjl} =  \mathbf{v}_{j,k}\mathbf{d}_{j,k}$, and  $\mathbf{v}_{j,l}\mathbf{d}_{jkjl} =  \mathbf{v}_{j,l}\mathbf{d}_{j,l}$,
so that $\mathbf{d}_{jkjl}$ is  in the intersection of  the directrix 
$ \mathbf{v}_{j,k}\mathbf{x} =  \mathbf{v}_{j,k}\mathbf{d}_{j,k}$ and the directrix  $\mathbf{v}_{j,l}\mathbf{x} =  \mathbf{v}_{j,l}\mathbf{d}_{j,l}$.  
Furthermore, $\mathbf{d}_{jkjl}$ is  in aff$(T)$, because it is a linear combination of $\mathbf{d}_{jk}$ and $\mathbf{u}_T$. \hfill $\Box$

\section{Intersecting the bisectors for all pairs of points in a subset of balls}

Given the set $\mathcal{P}$ of $m$  points in $\real^n$, and a Euclidean ball associated with each point in $\mathcal{P}$,   as defined in Section 9, let 
$\mathcal{S} =   \{\mathbf{p}_{i_1}, \ldots, \mathbf{p}_{i_s} \}$ be a subset of $\cal{P}$ whose centers are affinely independent.  Consider the problem of finding the intersection, denoted by $B_{\cal{S}}$, of bisectors $B_{i_j,i_k}$ for all pairs of points  $\mathbf{p}_{i_j}, \mathbf{p}_{i_k} \in \mathcal{S}$.  That is, find $B_{\mathcal{S}} = \cap_{1 \leq j < k \leq m} B_{i_j,i_k}$.  The intersection $B_{\mathcal{S}}$ is assumed to be non-empty.  If all the radii corresponding to points in $\mathcal{S}$ are equal,  $B_{\mathcal{S}}$ will be a hyperplane of dimension $n-s+1$,  and if some radii are unequal,  $B_{\mathcal{S}}$ will be a conic section of dimension $n-s+1$.
The solution approach uses the results of Section 9 to compute the parameters and vectors of $B_{\mathcal{S}}$ by intersecting one bisector with a sequence of hyperplanes.

The first property shows that $B_{\mathcal{S}}$ is equal to the intersection of only $s-1$ bisectors, $B_{i_1,i_j}$ for $j = 2,\ldots,s$.  The result holds for any other sequence of $s-1$ bisectors with the property that each point in $\mathcal{S}$ is in one of the $s-1$ pairs, and each pair is distinct.  The sequence above was chosen so that all the pairs have a common point, which simplifies the notation and exposition.

  \begin{property}
   $B_{\mathcal{S}} = \cap_{j=2}^{s} B_{i_1,i_j}$.
   \end{property}
   \noindent \textbf{Proof:}:
   Since $B_{\mathcal{S}}$ is the intersection of all bisectors corresponding to all pairs of points in $\mathcal{S}$, $B_{\mathcal{S}} \subseteq \cap_{j=2}^s B_{i_1,i_j}$. To show the opposite inclusion,
   observe that each point in $\mathcal{S}$ is associated with some bisector $B_{i_1,i_j}$ for   $j=2,\ldots,s$.
   Expression \eqref{bisec} implies 
   that $B_{i_1,i_j} \cap B_{i_1,i_k} \subseteq B_{i_j,i_k}$ for any triple $\{ \mathbf{p}_{i_1}, \mathbf{p}_{i_j}, \mathbf{p}_{i_k} \}$, 
   Then $\cap_{j=2}^{s} B_{i_1,i_j}
   = \cap_{2 \leq j < k \leq s} B_{i_1,i_j} \cap B_{i_1,i_k}   
   \subseteq \cap_{j=2}^s B_{i_1,i_j} \cap_{2 \leq j < k \leq s} B_{i_j,i_k} = B_{\mathcal{S}}$. 
    \hfill $\Box$ 
   
If the $s$ points in $\mathcal{S}$ have equal radii, then each bisector is a hyperplane and the intersection $B_{\mathcal{S}}$ is determined by solving the linear system corresponding to the hyperplanes $B_{i_1,i_j}$, given by \eqref{hyp}, for $j = 2,\ldots,s$.

If the points in $\mathcal{S}$ do not have equal radii, order the points in $\mathcal{S}$ by non-increasing radii, so that $\mathcal{S} = \{ \mathbf{p}_{i_1}, \ldots , \mathbf{p}_{i_{s}} \}$, with $r_{i_1} \geq  \ldots \geq r_{i_{s}}$.  By the assumption of unequal radii, $r_{i_1} > r_{i_{s}}$, so that the bisector $B_{i_1,i_s}$ is a hyperboloid.  In this case, $B_{\mathcal{S}}$ is constructed by intersecting $B_{i_1,i_s}$ with $s-1$ hyperplanes, which are constructed as follows.

For each of the $s-2$ triples of points $T_j = \{ \mathbf{p}_{i_1},\mathbf{p}_{i_j}, \mathbf{p}_{i_s} \}$, for $j = 2,\ldots,s-1$, Theorem 9.1 constructs the hyperplane $H_j$ and shows that $B_{i_1,i_s} \cap B_{i_1,i_j} = B_{i_1,i_s} \cap H_j$.
Thus $\cap_{j=2}^s B_{i_1,i_j} = B_{i_1,i_s} \cap_{j=2}^{s-1} H_j$, which 
is computed sequentially as $B_{i_1,i_s} \cap_{j=2}^k H_j$ for $k = 2,\ldots,s-1$. 
 
Initially, compute the vectors and parameters of $B_{i_1,i_s}$, and designate them as: $\mathbf{v}_{1}:=\mathbf{v}_{i_1,i_s}$, 
$\mathbf{c}_{1}:=\mathbf{c}_{i_1,i_s}$, $\mathbf{d}_{1}:=\mathbf{d}_{i_1,i_s}$, $\mathbf{a}_{1}:=\mathbf{a}_{i_1,i_s}$, 
$\epsilon_{1}:=\epsilon_{i_1,i_s}$, $a_{1}:=a_{i_1,i_s}$,  $b_{1}:=b_{i_1,i_s}$,  $c_{1}:=c_{i_1,i_s}$.
Also, define $\mathbf{hp}_1 = \mathbf{0}$.
Then for $k = 2,\ldots,s-1$, given the vectors and parameters $\mathbf{v}_{k-1}$, 
$\mathbf{c}_{k-1}$, $\mathbf{d}_{k-1}$, $\mathbf{a}_{k-1}$, 
$\epsilon_{k-1}$, $a_{k-1}$,  $b_{k-1}$,  $c_{k-1}$,
of $B_{i_1,i_s} \cap_{j=2}^{k-1} H_j$,
the expressions \eqref{HT12} through \eqref{pa} compute the vectors and parameters 
$\mathbf{v}_{k}$, 
$\mathbf{c}_{k}$, $\mathbf{d}_{k}$, $\mathbf{a}_{k}$, 
$\epsilon_{k}$, $a_{k}$,  $b_{k}$,  $c_{k}$,
of  $B_{i_1,i_s} \cap_{j=2}^{k-1}H_j \cap H_k$.
For each iteration $k = 3, \ldots, s-1$, $H_k$ must be projected onto $\cap_{j=2}^{k-1}H_j$.

\begin{alignat}{3}
\text{If\;} &r_{i_1} > r_{i_k} \text{\;\;\;\;\;\;} \mathbf{h}_{k} =   (\epsilon_{i_1,i_k}\mathbf{v}_{i_1,i_k} -  \epsilon_{i_1,i_{s}}\mathbf{v}_{i_1,i_{s}})/
\parallel \epsilon_{i_1,i_k}\mathbf{v}_{i_1,i_k} -  \epsilon_{i_1,i_{s}}\mathbf{v}_{i_1,i_{s}}\parallel    \label{HT12} \\
\text{If\;} &r_{i_1} = r_{i_k} \text{\;\;\;\;\;\;} \mathbf{h}_{k}  =  (\mathbf{p}_{i_1} -  \mathbf{p}_{i_k})/
\parallel \mathbf{p}_{i_1} -  \mathbf{p}_{i_k} \parallel   \label{HT2}\\
\mathbf{hp}_{k} &  = ( \mathbf{h}_{k} - \sum_{j=2}^{k-1}(\mathbf{h}_{k}  \mathbf{hp}_j) \mathbf{hp}_j )/
\parallel  \mathbf{h}_{k} - \sum_{j=2}^{k-1}(\mathbf{h}_{k}  \mathbf{hp}_j) \mathbf{hp}_j \parallel & \label{ph}\\
\mathbf{u}_{k-1} & = (\mathbf{hp}_{k} - (\mathbf{hp}_{k} \mathbf{v}_{k-1})\mathbf{v}_{k-1}) / 
\parallel \mathbf{hp}_{k} - (\mathbf{hp}_{k} \mathbf{v}_{k-1})\mathbf{v}_{k-1}\parallel   \label{pu}\\
\mathbf{v}_{k} & =  (\mathbf{v}_{k-1} - (\mathbf{v}_{k-1} \mathbf{hp}_{k})\mathbf{hp}_{k}) / 
\parallel \mathbf{v}_{k-1} - (\mathbf{v}_{k-1} \mathbf{hp}_{k})\mathbf{hp}_{k}  \parallel   \label{pvv}\\
\mathbf{d}_{k} & =  \mathbf{d}_{k-1} + \frac{\mathbf{v}_{i_1,i_k} (\mathbf{d}_{i_1,i_k} - \mathbf{d}_{k-1})}{\mathbf{v}_{i_1,i_k} \mathbf{u}_{k-1}} \mathbf{u}_{k-1}   \label{pd1}\\
\hat{h}_{k} & = \mathbf{hp}_{k} (\mathbf{d}_{k} - \mathbf{c}_{k-1} ) \label{phh} \\
\rho_{k} &= \mathbf{v}_{k-1}   \mathbf{v}_{k} \\
\sigma_{k} &= \mathbf{v}_{k-1}   \mathbf{hp}_{k} \\
\epsilon_{k} & =  \epsilon_{k-1}  \rho_{k} \\
\text{If\;}  &\epsilon_{k} \neq 1, \notag \\
& \tilde{c}_{k}  =  \epsilon_{k-1}^2 \rho_{k} \sigma_{k} \hat{h}_{k} / (1 - \epsilon_{k}^2) \\
& a_{k}^2 =  [(1 - \epsilon_{k-1}^2)(a_{k-1}^2(1 - \epsilon_{k}^2) - \hat{h}_{k}^2)] / (1 - \epsilon_{k}^2)^2 \label{lca} \\
& \mathbf{c}_{k}  =  \mathbf{c}_{k-1} + \hat{h}_{k}\mathbf{hp}_{k} + \tilde{c}_{k} \mathbf{v}_{k}  \label{pc}\\
& \mathbf{a}_{k}  =   \mathbf{c}_{k} + a_{k}\mathbf{v}_{k}   \label{pa} \\
\text{If\;}  &\epsilon_{k} > 1, \text{\;\;\;\;\;\;} b_{k}^2  =  - a_{k}^2(1 - \epsilon_{k}^2 ) = c_k^2 - a_k^2   \label{lcb}\\
\text{If\;} &\epsilon_{k} < 1, \text{\;\;\;\;\;\;} b_{k}^2  =   a_{k}^2(1 - \epsilon_{k}^2 )  = a_k^2-c_k^2  \label{lcba}\\
\text{If\;} &\epsilon_k = 1, \notag \\
&\tilde{c}_{k}  =  \epsilon_{k-1} \sigma_{k} \hat{h}_{k} / 2  \label{pch}\\
 &\hat{c}_{k} =  [(1 - \epsilon_{k-1}^2  \sigma_k^2) \hat{h}_{k}^2 + b_{k-1}^2 ] / (4 \tilde{c}_k) \label{pctld} \\
&\mathbf{c}_{k}  =  \mathbf{c}_{k-1} + \hat{h}_{k}\mathbf{hp}_{k} +\hat{c}_{k} \mathbf{v}_{k}  \label{pnc}
\end{alignat}
From expression (\ref{ph}), $\mathbf{hp}_2 = \mathbf{h}_2$,  and for each $k > 2$,  
$\mathbf{hp}_k$ is the orthogonal complement of the projection of $\mathbf{h}_{k}$ onto
the intersection of the hyperplanes $H_{2} \cap \ldots \cap H_{k-1}$.
Furthermore, the vectors $\mathbf{hp}_k$  are mutually orthogonal.
Given $\mathbf{hp}_k$, expression (\ref{pu}) computes the  vector $\mathbf{u}_{k-1}$ to be orthogonal to  $\mathbf{v}_{k-1}$,
and in the plane determined by  $\mathbf{hp}_k$ and  $\mathbf{v}_{k-1}$.
Expression (\ref{pvv}) computes the principal axis vector $\mathbf{v}_k$ for $B_{i_1,i_s} \cap_{j=2}^{k-1}H_j \cap H_k$.
Expression (\ref{pd1}) computes the point $\mathbf{d}_k$ which lies on the principal axis vector $\mathbf{v}_k$ 
through the center of $B_{i_1,i_s} \cap_{j=2}^{k-1}H_j \cap H_k$.
Expression (\ref{phh}) computes the distance from the center $\mathbf{c}_{k-1}$ of 
$B_{i_1,i_s} \cap_{j=2}^{k-1}H_j$ to the
point $\mathbf{d}_k$ on the principal axis of $B_{i_1,i_s} \cap_{j=2}^{k-1}H_j \cap H_k$.
If $\epsilon_{k-1}\rho_k \neq 1$, then the intersection is either a hyperboloid or an ellipsoid, and 
expression  (\ref{pc}) computes the center $\mathbf{c}_k$  based on   (\ref{intpar}),
and expression  (\ref{pa}) computes the vertices  $\mathbf{a}_k$ using  (\ref{intpar}).
If  $\epsilon_{k-1}\rho_k = 1$, then the intersection is a paraboloid, and 
expression  (\ref{pctld}) computes the center $\mathbf{c}_k$, which is also the vertex,  based on   (\ref{intpar}),
and expression  (\ref{pch}) computes the distance to the directrix or the focal point using  (\ref{intpar}).

The following property is useful in the application \cite{Dearingmincov}. 

\begin{property}\label{Aor}
For $k = 2, \dots, s-1$, $\mathbf{u}_k \mathbf{v}_1 = 0$, and $\mathbf{v}_k \mathbf{v}_1 >0$.
\end{property}
 \noindent \textbf{Proof:} Expression (\ref{ph}) implies that the vectors $\mathbf{hp}_k$, for $k = 2, \ldots, s-1$, are mutually orthogonal.
 Expression (\ref{pvv}) implies $\mathbf{v}_j  \mathbf{hp}_k = 0$, for $j \geq k, k=2, \ldots, s-1$.
 Expression  (\ref{pu}) implies $\mathbf{v}_{k-1} \mathbf{u}_{k-1} = 0$, for $ k=2, \ldots, s-1$,
 and  that $\mathbf{u}_{k}  \mathbf{hp}_{k} = 0$, for  $j \geq k, k=2, \ldots, s-1$.
 Expression (\ref{pu}) also implies  that 
 \begin{alignat}{1}
\mathbf{v}_{k} & =  \mathbf{vp}_{k} /\parallel  \mathbf{vp}_{k}  \parallel \text{ where    }   
 \mathbf{vp}_{k} =  \mathbf{v}_{1} - \sum_{j=2}^{k}(\mathbf{v}_{j-1} \mathbf{hp}_j) \mathbf{hp}_j. \label{pvp}
 \end{alignat}
Then $\mathbf{v}_{k}  \mathbf{u}_{k} = 0$ implies $\mathbf{vp}_{k}  \mathbf{u}_{k} = 0$, so that 
$0 =   \mathbf{u}_k\mathbf{v}_{1} - \sum_{j=2}^{k}(\mathbf{v}_{j-1}  \mathbf{hp}_j) \mathbf{u}_k \mathbf{hp}_j  \label{pv}$.
But $ \mathbf{u}_k\mathbf{h}_j = 0$ for $k \geq j$  so that $ \mathbf{u}_k\mathbf{v}_{0} = 0$ for $k= 2, \ldots, s-1$.

Finally, writing $\mathbf{v}_k \mathbf{v}_1$ using expression (\ref{pvp})  and combining with the results 
$\mathbf{v}_j  \mathbf{hp}_k = 0$, for $j \geq k,  k=2, \ldots, s-1$, gives
 $\mathbf{v}_k \mathbf{v}_1 =  \mathbf{v}_1 \mathbf{v}_1 = 1$,
 which gives the second conclusion. \;\;\;\;  \hfill $\Box$
 
 \section{Concluding comments}
 
The intersection of hyperboloids discussed here provides an alternative solution approach to the problem of Appolonius \cite{Brannan}.
Given three circles in a feasible configuration in $\real^2$, with centers $\mathbf{p}_1$, $\mathbf{p}_2$, $\mathbf{p}_3$, and positive radii $r_1, r_2, r_3$,  the problem of Appolonius is to construct eight circles, each of which is  internally or externally tangent to one, two or all three of the given circles in all eight combinations.  The vertex of the intersection of the bisectors of two of the pairs of given circles gives the center of the circle externally tangent to all three circles.  By multiplying the radii by $-1$, singularly, in pairs and all three, and adding a constant so that all resulting radii are positive,  and computing the intersection of the bisectors of two of the pairs gives each of the eight desired circles.  

This approach applies to the problem of Appolonius in higher dimensions as well.  

The quadratic form representations \eqref{HQ1}, \eqref{C} and \eqref{PQ} for hyperboloids and ellipsoids, cones, and paraboloids, respectively, has pedagogical advantages.  One  is that conic sections with any orientation (specified by the axis vector $\mathbf{v}$), center (specified by the point $\mathbf{c}$) and shape (specified by the eccentricity $\epsilon$) are represented by the quadratic form.  Another advantage, as shown by its use in this paper, is that the quadratic form is amenable to study and analysis of conic sections.

\end{document}